\documentclass[11pt, a4paper]{article}
\usepackage{latexsym,amssymb}
\usepackage{amsmath}
\usepackage[mathscr]{eucal}
\usepackage{color}
\usepackage[abbrev]{amsrefs}

\newtheorem{Theorem}{\bf Theorem}[section]
\newtheorem{Lemma}{\bf Lemma}[section]
\newtheorem{Proposition}{\bf Proposition}[section]
\newtheorem{Corollary}{\bf Corollary}[section]
\newtheorem{Remark}{\bf Remark}[section]
\newtheorem{Example}{\bf Example}[section]
\newtheorem{Definition}{\bf Definition}[section]

\newenvironment{theorem}{\begin{Theorem}$\!\!\!$}{\end{Theorem}}
\newenvironment{lemma}{\begin{Lemma}$\!\!\!$}{\end{Lemma}}
\newenvironment{proposition}{\begin{Proposition}$\!\!\!$}{\end{Proposition}}
\newenvironment{corollary}{\begin{Corollary}$\!\!\!$}{\end{Corollary}}
\newenvironment{remark}{\begin{Remark}$\!\!\!$}{\end{Remark}}

\newenvironment{definition}{\begin{Definition}$\!\!\!$}{\end{Definition}}

\numberwithin{equation}{section}

\numberwithin{equation}{section}

\begin{document}
\title{Eventual concavity properties of the heat flow}
\author{Kazuhiro Ishige\\
Graduate School of Mathematical Sciences, The University of Tokyo,\\ 
3-8-1 Komaba, Meguro-ku, Tokyo 153-8914, Japan\\
E-mail: {\tt ishige@ms.u-tokyo.ac.jp}}
\date{}
\maketitle
\begin{abstract}
The eventual concavity properties are useful to characterize geometric properties of the final state of solutions to parabolic equations. 
In this paper we give characterizations of the eventual concavity properties of the heat flow 
for nonnegative, bounded measurable initial functions with compact support.
\end{abstract}
\vspace{20pt}
\vspace{3pt}
\newpage
\section{Introduction}
Let $u$ be a nonnegative solution to the Cauchy problem for the heat equation 
$$
\partial_t u=\Delta u\quad\mbox{in}\quad{\mathbb R}^n\times(0,\infty),
\quad
u(\cdot,0)=\phi\ge 0\quad\mbox{in}\quad{\mathbb R}^n,
$$
where $\phi\in {\mathcal L}$ and $n\ge 1$.
Here and in what follows we denote by ${\mathcal L}$ the set of nontrivial, nonnegative, bounded measurable functions in ${\mathbb R}^n$ with compact support. 
The solution~$u$ is represented by 
\begin{equation}
\label{eq:1.1}
u(x,t)=(e^{t\Delta}\phi)(x)\equiv (4\pi t)^{-\frac{n}{2}}\int_{{\mathbb R}^n}e^{-\frac{|x-y|^2}{4t}}\phi(y)\,dy, 
\quad(x,t)\in{\mathbb R}^n\times(0,\infty),
\end{equation}
which leads to the eventual log-concavity property of $u$, 
that is, 
there exists $T>0$ such that 
$$
\log u((1-\mu)x+\mu y,t)\ge(1-\mu)\log u(x,t)+\mu \log u(y,t)
$$
for $x$, $y\in{\mathbb R}^n$, $\mu\in(0,1)$, and $t\in[T,\infty)$ 
(see e.g., \cite{LV}*{Theoerem~5.1}). 
It was proved in \cite{IST01}*{Theorem~3.3} that 
$e^{t\Delta}\phi$ eventually possesses stronger concavity properties than the log-concavity property (see also Theorem~\ref{Theorem:1.1}).
The eventual concavity properties are useful to characterize geometric properties of the final state 
of solutions to parabolic equations, and they have been studied for several nonlinear parabolic equations, 
e.g., a fully nonlinear variation of the heat equation (see \cite{HV}), 
the porous medium equation (see e.g., \cites{AV, BV, KL, LV}), and the evolution of $p$-Laplace equation (see \cite{LPV}). 

In this paper we characterize the eventual concavity properties of the function
\begin{equation*}
U_d[\phi](x,t):=(4\pi t)^d \left(e^{t\Delta}\phi\right)(x),
\quad (x,t)\in{\mathbb R}^n\times(0,\infty),
\end{equation*}
where $d\in{\mathbb R}$, using the notion of $F$-concavity (see Definition~\ref{Definition:1.1}), 
and study geometric properties of the final state of $U_d[\phi]$. 
$F$-concavity is a generalization of power concavity and, actually, the largest available generalization of the notion of concavity. 
Surprisingly, a variety of $F$-concavity properties which $U_d[\phi]$ eventually develops depends on $d$ and $n$. 
(See Theorem~\ref{Theorem:1.6} and Corollary~\ref{Corollary:1.3}.)
\subsection{$F$-concavity and the heat flow}
We follow \cite{IST05} to introduce the notion of $F$-concavity. 
\begin{definition}
 \label{Definition:1.1}
  Let $a\in(0,\infty]$, and set $I=[0,a)$ and $\mbox{{\rm int}}\,I=(0,a)$. 
 \begin{itemize}
  \item[{\rm (1)}] 
  A function~$F:I\to[-\infty,\infty)$ is said admissible on $I$ if $F\in C(\mbox{{\rm int}}\,I)$, 
  $F$ is strictly increasing on $I$, and $F(0)=-\infty$. 
  Let $\mathcal {A}(I)$ be the set of $I$-valued functions in ${\mathbb R}^n$.  
  \item[{\rm (2)}] 
  Let $F$ be admissible on $I$.  
  For any $f\in{\mathcal A}(I)$, we say that $f$ is $F$-concave in ${\mathbb R}^n$ 
  if 
  $$
  F(f((1-\mu)x+\mu y))\ge(1-\mu)F(f(x))+\mu F(f(y))
  $$
  for $x,y\in{\mathbb R}^n$ and $\mu\in (0,1)$. 
  We denote by ${\mathcal C}[F]$ the set of $F$-concave functions in ${\mathbb R}^n$. 
  \item[{\rm (3)}]
  Let $F_1$ and $F_2$ be admissible on $I$. 
  We say that $F_1$-concavity is weaker {\rm({\it resp.}\,\,strictly weaker)} than $F_2$-concavity in ${\mathcal A}(I)$, 
  or equivalently that  $F_2$-concavity is stronger {\rm({\it resp.}\,\,strictly stronger)} than $F_1$-concavity in ${\mathcal A}(I)$ if 
  $\mathcal{C}[F_2]\subset\mathcal{C}[F_1]$ 
  {\rm (}resp.\,\,$\mathcal{C}[F_2]\subsetneq\mathcal{C}[F_1]${\rm)}.
\end{itemize}
\end{definition}
We define the eventual $F$-concavity property of $U_d[\phi]$ and 
the preservation of the $F$-concavity property by the heat flow in ${\mathbb R}^n$.
 \begin{definition}
 \label{Definition:1.2}
 Let $F$ be admissible on $I=[0,a)$ with $a\in(0,\infty]$. 
  \begin{itemize}
  \item[{\rm (1)}] 
  Let $d\in{\mathbb R}$ and $\phi\in{\mathcal L}$.   
  We say that $U_d[\phi]$ possesses the eventual $F$-concavity property if 
  there exists $T>0$ such that $U_d[\phi](\cdot,t)$ is $F$-concave in ${\mathbb R}^n$ for $t\in[T,\infty)$. 
  \item[{\rm (2)}] 
  We say that the $F$-concavity property is preserved by the heat flow in ${\mathbb R}^n$ if, 
  for any $\phi\in{\mathcal C}[F]\cap L^\infty({\mathbb R}^n)$, 
  $e^{t\Delta}\phi$ is $F$-concave in ${\mathbb R}^n$ for $t>0$.
\end{itemize}
 \end{definition}
The preservation of concavity properties in parabolic problems 
has been studied in many papers since the pioneering work by Brascamp--Lieb~\cite{BL} 
(see e.g., \cites{BV, CW1, CW2, CW3, GK, INS, IS01, IS02, IST01, IST03, IST05, LV, Keady, Ken, Ken02, Kol, Kor} and references therein). 
Among others, in \cites{IST05}, 
the author of this paper, Salani, and Takatsu gave a complete characterization of the $F$-concavity property preserved by 
the heat flow in ${\mathbb R}^n$, more generally, by the Dirichlet heat flow in convex domains in ${\mathbb R}^n$. 
The preservation of the $F$-concavity property by the heat flow causes the following question. 
\begin{itemize}
  \item[(Q)] 
  Assume that the $F$-concavity property is preserved by the heat flow in ${\mathbb R}^n$. 
  Then does $e^{t\Delta}\phi$ possess the eventual $F$-concavity for $\phi\in{\mathcal L}$ 
  without the $F$-concavity property of the initial data~$\phi$?  
\end{itemize}
One of the aims of this paper is to give an affirmative answer to question~(Q) 
for the initial data $\phi$ in a subclass ${\mathcal L}_A$ (see Section~1.3) of ${\mathcal L}$. 

We give some typical examples of $F$-concavity, and recall some results on the preservation of concavity properties by the heat flow in ${\mathbb R}^n$.
\begin{itemize}
  \item[(E1)] (Power concavity)\\
  Let $I=[0,\infty)$ and $\alpha\in{\mathbb R}$. Let $\Phi_\alpha$ be an admissible function on $I$ defined by 
  \begin{equation*}
  \Phi_\alpha(\tau):=
  \left\{
  \begin{array}{ll}
  \displaystyle{\frac{\tau^\alpha-1}{\alpha}}\quad & \mbox{for $\tau>0$ if $\alpha\not=1$},\vspace{7pt}\\
  \log\tau\quad & \mbox{for $\tau>0$ if $\alpha=1$},\vspace{5pt}\\
 -\infty\quad & \mbox{for $\tau=0$}.
  \end{array}
  \right.
  \end{equation*}
  A nonnegative function $f$ in ${\mathbb R}^n$ is said $\alpha$-concave in ${\mathbb R}^n$ 
  if $f$ is $\Phi_\alpha$-concave in~${\mathbb R}^n$.
  We often call $0$-concavity log-concavity. 
  Power concavity is a generic term of $\alpha$-concavity, and it possesses the following properties.
  \begin{itemize}
  \item[(a)] 
  If $f$ is $\alpha$-concave in ${\mathbb R}^n$ and $\alpha'\le\alpha$, then $f$ is $\alpha'$-concave in ${\mathbb R}^n$.
  \item[(b)] 
  If $f$ is $\alpha$-concave in ${\mathbb R}^n$ and $\lambda\in(0,\infty)$, then $\lambda f$ is $\alpha$-concave in ${\mathbb R}^n$.
  \item[(c)]
  Among power concavity properties, only the log-concavity property is preserved by the heat flow in ${\mathbb R}^n$.
  \end{itemize}
  Properties~(a) and (c) follow from the Jensen inequality and \cite{IST05}*{Theorem~1.5}, respectively. 
  Property~(b) is a characteristic property of power concavity 
  (see Lemma~\ref{Lemma:2.3}).  
  \item[(E2)] (Power log-concavity)\\
  Let $I=[0,1)$ and $\beta\in(0,\infty)$. 
  Let $\Psi_\beta$ be an admissible function on $I$ defined by 
  \begin{equation}
  \label{eq:1.2}
  \Psi_\beta(\tau)\equiv 
  \left\{
  \begin{array}{ll}
  -\left(-\log\tau\right)^\beta & \mbox{for $\tau\in(0,1)$},\vspace{3pt}\\
  -\infty & \mbox{for $\tau=0$}.
  \end{array}
  \right.
  \end{equation}
  A function $f\in{\mathcal A}(I)$ is said $\beta$-log-concave in ${\mathbb R}^n$
  if $f$ is $\Psi_\beta$-concave in~${\mathbb R}^n$. 
  Power concavity is a generic term of $\beta$-log-concavity, and 
  it possesses the following properties (see \cites{IST01, IST03}).
  \begin{itemize}
  \item[(a)] 
  If $f$ is $\beta$-log-concave in ${\mathbb R}^n$ and $\beta\le\beta'$, then $f$ is $\beta'$-log-concave in ${\mathbb R}^n$.
  \item[(b)] 
  If $f$ is $\beta$-log-concave in ${\mathbb R}^n$ with $\beta\in(0,1]$ and $\lambda\in(0,1]$, then $\lambda f$ is $\beta$-log-concave in~${\mathbb R}^n$.
  \item[(c)] 
  The $\beta$-log-concavity property is preserved by the heat flow  if and only if $\beta\in[1/2,1]$
  (see \cite{IST01}*{Theorem~3.1} and \cite{IST05}*{Theorem~6.9}).
  \end{itemize}
  \item[(E3)] (Hot-concavity)\\
  Let $h$ be a positive function in ${\mathbb R}$ defined by
  $$
  h(x):=(4\pi)^{-\frac{1}{2}}\int_0^\infty e^{-\frac{(x-z)^2}{4}}\,dz,\quad x\in{\mathbb R}.
  $$
  Then $h$ is strictly increasing in ${\mathbb R}$, $\lim_{x\to-\infty}h(x)=0$, and $\lim_{x\to\infty}h(x)=1$. 
  Let $H$ be the inverse function of $h$, which is defined in $(0,1)$, and set $H(0)=-\infty$. 
  Then $H$ is admissible on $I=[0,1)$, and $H$-concavity is called hot-concavity
  (see \cites{IST05, IST06}).   
  Hot-concavity possesses the following properties (see \cite{IST05}*{Theorems~5.1 and 5.2}). 
  \begin{itemize}
  \item[(a)]
  If $f$ is hot-concave in ${\mathbb R}^n$ and $\lambda\in(0,1]$, then $\lambda f$ is hot-concave in ${\mathbb R}^n$.
  \item[(b)] 
  The hot-concavity property is the strongest $F$-concavity properties in ${\mathcal A}(I)$ preserved by the heat flow in ${\mathbb R}^n$.
  \end{itemize}
\end{itemize}
Applying \cite{IST01}*{Theorem~3.3} with (E2)-(a), 
we see that the heat flow possesses the following eventual concavity properties. 
\begin{theorem}
\label{Theorem:1.1}
Let $\phi\in{\mathcal L}$. 
For any $\beta>1/2$, $e^{t\Delta}\phi$ possesses the eventual $\beta$-log-concavity concavity, 
that is, 
there exists $T>0$ such that 
$$
\Psi_\beta\left((e^{t\Delta}\phi)((1-\mu)x+\lambda y)\right)\ge(1-\mu)\Psi_\beta((e^{t\Delta}\phi)(x))+\mu\Psi_\beta((e^{t\Delta}\phi)(y))
$$
for $x$, $y\in{\mathbb R}^n$, $\mu\in(0,1)$, and $t\in[T,\infty)$, where $\Psi_\beta$ is as in \eqref{eq:1.2}.
\end{theorem}
While Theorem~\ref{Theorem:1.1} give a partial affirmative answer to question~(Q), 
the following question~(Q') naturally arises.
\begin{itemize}
  \item[(Q')] 
  What kind of eventual $F$-concavity properties does the heat flow develops for $\phi\in{\mathcal L}$? 
  In particular, does $e^{t\Delta}\phi$ possess the eventual $1/2$-log-concavity concavity property and the eventual hot-concavity 
  for $\phi\in{\mathcal L}$ ? 
\end{itemize}
In this paper 
we investigate the eventual $F$-concavity properties of $U_d[\phi]$ 
and reveal the relation between the eventual $F$-concavity properties of $U_d[\phi]$, the behavior of $F(\tau)$ as $\tau\to+0$, and the parameter~$d$.
Furthermore, as an application of our study, we solve questions~(Q) and (Q') for initial data $\phi$ in a suitable subclass ${\mathcal L}_A$ of ${\mathcal L}$. 
\subsection{Main results in the class ${\mathcal L}$}
In this section we state our results on the eventual $F$-concavity properties of $U_d[\phi]$ for $\phi\in{\mathcal L}$. 
In the case of $d<n/2$, we assume the following condition 
with $a\in(0,\infty)$. 
\begin{itemize}
  \item[($\mbox{F}_a$)] 
  Let $F$ be admissible on $I=[0,a)$ such that $F\in C^2(\mbox{int}\,I)$ and $F'>0$ in $\mbox{int}\,I$. 
\end{itemize}
In the case of $d<n/2$, 
since $\|U_d(t)\|_{L^\infty({\mathbb R}^n)}\to 0$ as $t\to\infty$, 
it is not restrictive to assume that $a<\infty$. 
Set
\begin{equation}
\label{eq:1.3}
{\mathcal F}(r)=F(ae^{-r})\quad\mbox{and}\quad
\kappa_F(r):=-r\frac{{\mathcal F}''(r)}{{\mathcal F}'(r)}\quad\mbox{for $r\in(0,\infty)$}. 
\end{equation}
The function~$\kappa_F$ is independent of $a\in(0,\infty)$ (see Remark~\ref{Remark:1.1}~(1)), and it has the following property: 
when $F_1$ and $F_2$ satisfies condition~($\mbox{F}_a$) with $a\in(0,\infty)$, 
\begin{equation}
\label{eq:1.4}
\mbox{$\kappa_{F_1}\ge\kappa_{F_2}$ in $(0,\infty)$ if and only if $F_1$-concavity is stronger than $F_2$-concavity in ${\mathcal A}(I)$}.
\end{equation}
(See Lemma~\ref{Lemma:2.6}.) 
The first result of this paper concerns with a characterization of the eventual $F$-concavity of $U_d[\phi]$ with $d<n/2$. 
\begin{theorem}
\label{Theorem:1.2}
Let $d<n/2$. 
Assume condition~{\rm ($\mbox{F}_a$)} with $a\in(0,\infty)$. Set 
$$
\kappa_F^*:=\displaystyle{\limsup_{r\to\infty}}\,\kappa_F(r).
$$ 
\begin{itemize}
  \item[{\rm (1)}] 
  If $\kappa_F^*<1/2$, then 
  $U_d[\phi]$  possesses the eventual $F$-concavity property for all $\phi\in{\mathcal L}$. 
  \item[{\rm (2)}] 
  If $\kappa_F^*>1/2$, then
  $U_d[\phi]$ does not possess the eventual $F$-concavity property for any $\phi\in{\mathcal L}$. 
\end{itemize}
\end{theorem}
The case of $\kappa_F^*=1/2$ is delicate, and it is treated in Section 1.3. 
\begin{remark}
\label{Remark:1.1}
We give some comments on the function~$\kappa_F$.
\begin{itemize}
  \item[{\rm (1)}] 
  Let $F$ satisfy condition~{\rm ($\mbox{F}_a$)} with $a\in(0,\infty)$. For any $k>0$, set $F_k(r):=F(kr)$ for $r\in [0,k^{-1}a)$. 
  Then $F_k$ is admissible on $[0,k^{-1}a)$ and 
  $$
  {\mathcal F}_k(r):=F_k(k^{-1}ae^{-r})=F(ae^{-r})={\mathcal F}(r),\quad 
  \kappa_{F_k}(r)=\kappa_F(r),
  $$
  for $r\in(0,\infty)$. 
  In particular, $\kappa_F^*=\kappa_{F_k}^*$ for $k>0$.  
  \item[{\rm (2)}] {\rm({\it Power concavity})}\\
  Let $F=\Phi_\alpha(\tau)$ for $\tau\in[0,1)$, where $\alpha\in{\mathbb R}$ and $\Phi_\alpha$ is as in {\rm (E1)}.
  Since
  $$
  {\mathcal F}(r)=F(e^{-r})=
  \left\{
  \begin{array}{ll}
  \displaystyle{\frac{e^{-\alpha r}-1}{\alpha}} & \mbox{if}\quad\alpha\not=0,\vspace{5pt}\\
  -r & \mbox{if}\quad\alpha=0,
  \end{array}
  \right.
  \quad
  \kappa_F(r)=\alpha r,\quad\mbox{for $r\in(0,\infty)$},
  $$
  we see that $\kappa_F^*=\infty$ if $\alpha>0$, $\kappa_F^*=0$ if $\alpha=0$, and $\kappa_F^*=-\infty$ if $\alpha<0$. 
  Then it follows from Theorem~{\rm\ref{Theorem:1.2}} that, for any $\phi\in{\mathcal L}$, 
  $U_d[\phi]$ possesses the eventual $\alpha$-concavity property if and only if $\alpha\le 0$.
  \item[{\rm (3)}] {\rm({\it Power log-concavity})}\\
  Let $F=\Psi_\beta(\tau)$ for $\tau\in[0,1)$, where $\beta\in(0,\infty)$ and $\Psi_\beta$ is as in {\rm (E2)}.
  Then
  $$
  {\mathcal F}(r)=F(e^{-r})=-r^\beta,
  \quad
  \kappa_F(r)=-\beta+1,\quad\mbox{for $r\in(0,\infty)$}. 
  $$
  This implies that $\kappa^*_F<1/2$ if and only if $\beta>1/2$.
  \item[{\rm (4)}]  
  Let $F$ satisfy condition~{\rm ($\mbox{F}_a$)} with $a=1$. 
  Then $F$-concavity is stronger than $\beta$-log-concavity in ${\mathcal A}([0,1))$ if and only if 
  $$
  \kappa_F(r)\ge -\beta+1,\quad r\in(0,\infty).
  $$
  This follows from \eqref{eq:1.4} and Remark~{\rm \ref{Remark:1.1}}~{\rm (3)}.
  \item[{\rm (5)}] 
  Let $F$ satisfy condition~{\rm ($\mbox{F}_a$)} with $a\in(0,\infty)$. 
  If the $F$-concavity property is preserved by the heat flow, then $\kappa_F^*\le 1/2$ {\rm ({\it see} Theorem~\ref{Theorem:1.6}~(1))}. 
  On the other hand, 
  the hot-concavity property is preserved by the heat flow, 
  and hot-concavity is stronger than $1/2$-log-concavity in ${\mathcal A}([0,1))$ {\rm ({\it see} (E2)-(c) {\it and} (E3)-(b))}. 
  Then we observe from Remark~{\rm \ref{Remark:1.1}}~{\rm (4)} that
  $$
  \kappa_F^*=\frac{1}{2}\quad\mbox{in the case of hot-concavity}. 
  $$
  \end{itemize}
\end{remark}

Next, we consider the case of $d=n/2$. 
Let $F$ be admissible on $I=[0,a)$ with $a\in(0,\infty]$.  
Let $\phi\in{\mathcal L}$ be such that 
$$
M_\phi:=\int_{{\mathbb R}^n}\phi(y)\,dy<a. 
$$
Then
\begin{equation}
\label{eq:1.5}
U_{\frac{n}{2}}[\phi](x,t)=\int_{{\mathbb R}^n}e^{-\frac{|x-y|^2}{4t}}\phi(y)\,dy\le M_\phi<a,
\quad (x,t)\in{\mathbb R}^n\times(0,\infty),
\end{equation}
and $F(U_d[\phi])$ is well-defined in ${\mathbb R}^n\times(0,\infty)$. 
%
\begin{theorem}
\label{Theorem:1.3}
Let $d=n/2$. Assume condition~{\rm ($\mbox{F}_a$)} with $a\in(0,\infty)$.
\begin{itemize}
  \item[{\rm (1)}] 
  Assume that 
  \begin{equation}
  \label{eq:1.6}
  \kappa_F(r)\le\frac{1}{2}\quad\mbox{for $r\in(0,\infty)$},
  \qquad
  \kappa_F^*<\frac{1}{2}.
  \end{equation}
  Then $U_d[\phi]$  possesses the eventual $F$-concavity property for all $\phi\in{\mathcal L}$ with $M_\phi<a$. 
  \item[{\rm (2)}] 
  Assume that $\kappa_F(r_*)>1/2$
  for some $r_*\in(0,\infty)$. 
  Then $U_d[\phi]$ does not possess the eventual $F$-concavity property for some $\phi\in{\mathcal L}$ with $M_\phi<a$. 
\end{itemize}
\end{theorem}
\begin{theorem}
\label{Theorem:1.4}
Let $F$ be admissible on $[0,\infty)$ and $d=n/2$. 
Then $U_d[\phi]$  possesses the eventual $F$-concavity property for all $\phi\in{\mathcal L}$ if and only if $F$-concavity is weaker than log-concavity in ${\mathcal A}([0,\infty))$. 
\end{theorem}
Applying Theorems~\ref{Theorem:1.2} and \ref{Theorem:1.3} to power log-concavity, we have:
\begin{corollary}
\label{Corollary:1.1}
\begin{itemize}
  \item[{\rm (1)}] 
  Let $d<n/2$. 
  \begin{itemize}
  \item[{\rm (a)}]
  If $\beta>1/2$, then $U_d[\phi]$ possesses the eventual $\beta$-log-concavity property for all $\phi\in{\mathcal L}$. 
  \item[{\rm (b)}]
  If $\beta<1/2$, then $U_d[\phi]$ does not possess the eventual $\beta$-log-concavity property for any $\phi\in{\mathcal L}$.
  \end{itemize}
  \item[{\rm (2)}] 
  Let $d=n/2$. 
  \begin{itemize}
  \item[{\rm (a)}]
  If $\beta>1/2$, then $U_d[\phi]$ possesses the eventual $\beta$-log-concavity property for all $\phi\in{\mathcal L}$ with $M_\phi<1$. 
  \item[{\rm (b)}]
  If $\beta<1/2$, then $U_d[\phi]$ does not possess the eventual $\beta$-log-concavity property for some $\phi\in{\mathcal L}$ with $M_\phi<1$.
  \end{itemize}
\end{itemize}
\end{corollary}

Next, we consider the case of $d>n/2$. 
Then $\|U_d(t)\|_{L^\infty({\mathbb R}^n)}\to\infty$ as $t\to\infty$, 
which implies that $F(U_d[\phi](\cdot,t))$ is well-defined for all large enough $t$  if and only if $a=\infty$.
\begin{theorem}
\label{Theorem:1.5}
Let $F$ be admissible on $[0,\infty)$, $d>n/2$, and $\phi\in {\mathcal L}$. 
Then 
$U_d[\phi]$ possesses the eventual $F$-concavity property if and only if $F$-concavity is weaker than log-concavity in ${\mathcal A}([0,\infty))$. 
\end{theorem}
Theorems~\ref{Theorem:1.4} and \ref{Theorem:1.5} completely characterize the eventual $F$-concavity properties
which $U_d[\phi]$ develops for all $\phi\in{\mathcal L}$ in the case of $a=\infty$ and $d\ge n/2$. 
In the next section we introduce a subclass ${\mathcal L}_A$ of ${\mathcal L}$, 
and study the eventual $F$-concavity properties which $U_d[\phi]$ develops for all $\phi\in{\mathcal L}_A$ in the case of $a\in(0,\infty)$ and $d\le n/2$.
\subsection{Main results in the class ${\mathcal L}_A$}
For any $\phi\in{\mathcal L}$, 
we say that $\phi\in{\mathcal L}_A$ if $\phi$ satisfies the following condition.
\begin{itemize}
  \item[(A)] 
  There exists $C>0$ such that
  $$
  \int_{{\mathbb R}^n}\langle y,\xi\rangle^2e^{k\langle y,\xi\rangle}\phi(y-z)\,dy
  \le Ck^{-2}\int_{{\mathbb R}^n}e^{k\langle y,\xi\rangle}\phi(y-z)\,dy
  $$
  for all $k\ge 1$ and $(\xi,z)\in{\mathbb S}^{n-1}\times{\mathbb R}^n$ satisfying $\mbox{ess sup}_{y\in P_\phi}\langle y-z,\xi\rangle=0$.  
  Here $P_\phi$ is the positive set of $\phi$, that is, 
  $P_\phi:=\{x\in{\mathbb R}^n\,:\,\phi(x)>0\}$. 
\end{itemize} 
See Section~2.3 for sufficient conditions for $\phi\in{\mathcal L}$ to satisfy condition~(A). 
For example, if there exist $C>0$ and a bounded smooth open set $\Omega$ in ${\mathbb R}^n$ such that 
$$
C^{-1}\chi_\Omega(x)\le\phi(x)\le C\chi_\Omega(x)\quad\mbox{for a.a.~$x\in{\mathbb R}^n$},
$$
then $\phi\in{\mathcal L}_A$. Here $\chi_\Omega$ is the characteristic function of $\Omega$.

In this subsection we give a complete characterization of the eventual $F$-concavity properties which $U_d[\phi]$ develops for all $\phi\in{\mathcal L}_A$
in the case of $d\le n/2$. The behavior of the function $\sigma_F$ defined by 
\begin{equation}
\label{eq:1.7}
\sigma_F(r):=r\left(\kappa_F(r)-\frac{1}{2}\right),\quad r\in(0,\infty),
\end{equation}
is crucial for the characterization.
\begin{theorem}
\label{Theorem:1.6}
Assume condition~{\rm ($\mbox{F}_a$)} with $a\in(0,\infty)$. 
Let $\sigma_F$ be as in \eqref{eq:1.7}. 
\begin{itemize}
  \item[{\rm (1)}] 
  If the $F$-concavity property is preserved by the heat flow in ${\mathbb R}^n$, then 
  $\displaystyle{\limsup_{r\to\infty}}\,\sigma_F(r)<\infty$ and $\kappa_F^*\le 1/2$. 
  \item[{\rm (2)}] 
  Let $d<(n-1)/2$. 
  Then $U_d[\phi]$ possesses the eventual $F$-concavity property for all $\phi\in{\mathcal L}_A$ if and only if 
  \begin{equation}
  \label{eq:1.8}
  \limsup_{r\to\infty}\left(\sigma_F(r)-\frac{1}{4}\log r\right)<\infty.
  \end{equation}
  \item[{\rm (3)}] 
  Let $(n-1)/2\le d<n/2$. 
  Then $U_d[\phi]$ possesses the eventual $F$-concavity property for all $\phi\in{\mathcal L}_A$ if and only if 
  \begin{equation}
  \label{eq:1.9}
  \lim_{r\to\infty}\left(\sigma_F(r)-\frac{1}{2}\left(\frac{n}{2}-d\right)\log r\right)=-\infty.
  \end{equation}
  \item[{\rm (4)}] 
  Let $d=n/2$ and $M_\phi<a$. 
  Then $U_d[\phi]$ possesses the eventual $F$-concavity property for all $\phi\in{\mathcal L}_A$ if and only if 
  \begin{equation}
  \label{eq:1.10}
  \kappa_F(r)\le\frac{1}{2}\quad\mbox{for $r\in(0,\infty)$},
  \qquad
  \lim_{r\to\infty}\sigma_F(r)=-\infty.
  \end{equation}
 \end{itemize}
\end{theorem}
We remark that 
\begin{align*}
\mbox{\eqref{eq:1.10} holds}\qquad & \Rightarrow\qquad\mbox{\eqref{eq:1.9} holds with $(n-1)\le d<n/2$}\\
& \Rightarrow\qquad\mbox{\eqref{eq:1.8} holds with $d<(n-1)/2$}.
\end{align*}
Furthermore, by Theorems~\ref{Theorem:1.2}--\ref{Theorem:1.6} 
we have the following list on whether $U_d[\phi]$ possesses the eventual $F$-concavity property for all $\phi\in{\mathcal L}_A$ or not. 

$$
\begin{tabular}{|c||c|c|c|}
 \hline
  & $\kappa_F^*<1/2$ & $\kappa_F^*=1/2$ & $\kappa_F^*>1/2$\\
 \hline \hline
$2d<n-1$ & Yes  & Yes iff \eqref{eq:1.8} holds & No \\
 \hline
 $n-1\le 2d<n$
 & Yes & Yes iff \eqref{eq:1.9} holds & No \\
 \hline
 $2d=n,\,M_\phi<a\in(0,\infty)$ & Yes iff $\kappa_F(r)\le 1/2$  & Yes iff \eqref{eq:1.10} holds & No \\
  \hline\hline
 $2d\ge n,\,a=\infty$ & 
 \multicolumn{3}{|c|}{Yes iff $F$-concavity is weaker than log-concavity in ${\mathcal A}([0,\infty))$}\\
  \hline
\end{tabular}
$$
\vspace{3pt}

\noindent
As a direct consequence of Theorem~\ref{Theorem:1.6}, 
we have the following results, which give our answers to questions~(Q) and (Q') for $\phi\in{\mathcal L}_A$. 
\begin{corollary}
\label{Corollary:1.2}
\begin{itemize}
 \item[{\rm (1)}]
  Let $F$ be admissible on $I=[0,a)$ with $a\in(0,\infty]$ and $d<n/2$.
  If the $F$-concavity property is preserved by the heat flow in ${\mathbb R}^n$, then $U_d[\phi]$ possesses the eventual $F$-concavity property
  for all $\phi\in{\mathcal L}_A$. 
  \item[{\rm (2)}]
  Let $d<n/2$. 
  Then $U_d[\phi]$ possesses the eventual $1/2$-log-concavity and the eventual hot-concavity property for all $\phi\in{\mathcal L}_A$.
  \item[{\rm (3)}]
  Let $d=n/2$. 
  Then $U_d[\phi]$ does not possess the eventual $1/2$-log-concavity property and the eventual hot-concavity property for some $\phi\in{\mathcal L}_A$ with $M_\phi<1$.
\end{itemize}
\end{corollary}
\begin{corollary}
\label{Corollary:1.3} 
Assume condition~{\rm ($\mbox{F}_a$)} with $a\in(0,\infty)$.
Then the heat fow
$e^{t\Delta}\phi$ possesses the eventual $F$-concavity property for all $\phi\in{\mathcal L}_A$ 
if and only if 
\begin{equation*}
\left\{
\begin{array}{ll}
\displaystyle{\limsup_{r\to\infty}\left(\sigma_F(r)-\frac{1}{4}\log r\right)<\infty} & \quad\mbox{in the case of $n\ge 2$},\vspace{7pt}\\
\displaystyle{\lim_{r\to\infty}\left(\sigma_F(r)-\frac{1}{4}\log r\right)=-\infty} & \quad\mbox{in the case of $n=1$}\vspace{3pt}.
\end{array}
\right.
\end{equation*}
In particular, there exists an admissible function $F$ in $I=[0,1)$, satisfying condition~{\rm ($\mbox{F}_a$)} with $a=1$, such that 
\begin{itemize}
  \item[{\rm (1)}] 
  $F$-concavity is strictly stronger than hot-concavity in ${\mathcal A}([0,1))$;
  \item[{\rm (2)}] 
  the $F$-concavity property is not preserved by the heat flow in ${\mathbb R}^n$;
  \item[{\rm (3)}] 
  $e^{t\Delta}\phi$ possesses the eventual $F$-concavity property for all $\phi\in{\mathcal L}_A$. 
\end{itemize}
\end{corollary}
By Corollaries~\ref{Corollary:1.1} and \ref{Corollary:1.2} 
we have the following list on whether $U_d[\phi]$ possesses the eventual $\beta$-concavity property for all $\phi\in{\mathcal L}_A$ or not.

$$
\begin{tabular}{|c||c|c|c|}
 \hline
$\beta$-log-concavity  & $\beta>1/2$ & $\beta=1/2$ & $\beta<1/2$\\
 \hline \hline
$2d<n$ & Yes  & Yes  & No \\
 \hline
 $2d=n,\,M_\phi<1$ & Yes  & No & No \\
  \hline
\end{tabular}
$$
\vspace{3pt}

We explain our strategy of the proofs of our main results. 
In our analysis it is crucial to study geometric properties of $U_d[\phi]$ outside parabolic cones $B(0,L\sqrt{t})$, 
in particular, outside $B(0,Lt)$, for large enough $t$, where $L>0$.  
Let $\phi\in{\mathcal L}$. 
In the proofs of Theorem~\ref{Theorem:1.2}, Theorem~\ref{Theorem:1.3}, and Theorem~\ref{Theorem:1.6}~(2)--(4), 
we discuss the eventual $F$-concavity property of $U_d[\phi]$ studying the sign of the function
$$
\frac{\partial^2}{\partial\xi^2}F(U_d[\phi](x,t))\equiv \frac{\partial^2}{\partial\xi^2}{\mathcal F}(v_d(x,t)),
\quad \xi\in{\mathbb S}^{n-1},
$$
for $x\in{\mathbb R}^n$ and large enough $t$, where
\begin{align}
\label{eq:1.11}
v_d(x,t) & :=-\log \left(U_d[a^{-1}\phi](x,t)\right)=\left(\frac{n}{2}-d\right)\log(4\pi t)+\frac{|x|^2}{4t}-\log w(x,t),\\
\label{eq:1.12}
w(x,t) & :=\int_{{\mathbb R}^n}e^{\frac{2\langle x,y\rangle-|y|^2}{4t}}\,a^{-1}\phi(y)\,dy.
\end{align}
It easily follows from $\phi\in{\mathcal L}$ that 
\begin{equation}
\label{eq:1.13a}
\lim_{|x|\to\infty}v_d(x,t)=\infty\quad\mbox{for $t\in(0,\infty)$}
\end{equation}
and
\begin{equation}
\label{eq:1.13}
\left|\left(\frac{\partial}{\partial\xi}\log w\right)(x,t)\right|\le Ct^{-1},
\quad
\left|\left(\frac{\partial^2}{\partial\xi^2}\log w\right)(x,t)\right|\le Ct^{-2},
\end{equation}
for $(x,t)\in{\mathbb R}^n\times(0,\infty)$ and $\xi\in{\mathbb S}^{n-1}$
(see Lemma~\ref{Lemma:2.8}), where $C>0$.
In the case of $\kappa_F^*\not=1/2$, thanks to the function $\kappa_F$ and the term $|x|^2/4t$ in \eqref{eq:1.11}, we study the sign of 
$\partial^2 {\mathcal F}(v_d)/\partial\xi^2$
systematically (see Proposition~\ref{Proposition:3.1}), and prove Theorems~\ref{Theorem:1.2} and \ref{Theorem:1.3}.
On the other hand, in the case of $\kappa_F^*=1/2$, 
we require more precise decay estimates of the derivatives of $w$ in ${\mathbb R}^n\setminus B(0,Lt)$ for large enough $t$, where $L>0$, than those of \eqref{eq:1.13}.
Indeed, we improve decay estimate~\eqref{eq:1.13} under condition~(A) (see Lemma~\ref{Lemma:2.9}). 
Furthermore, introducing the functions $\nu_F$ and $\psi_d$ (see \eqref{eq:2.3} and \eqref{eq:4.9}, respectively) and modifying the proof in the case of $\kappa_F^*\not=1/2$,  
we study the sign of $\partial^2 {\mathcal F}(v_d)/\partial\xi^2$ in the case of $\kappa_F^*=1/2$ (see Proposition~\ref{Proposition:4.1}) to prove Theorem~\ref{Theorem:1.6}. 
We also develop the arguments in \cites{IST03} to 
prove Theorems~\ref{Theorem:1.4} and \ref{Theorem:1.5}.
\vspace{5pt}

The rest of this paper is organized as follows. 
Section~2 is divided into three subsections. 
In Section~2.1 we recall some properties of $F$-concavity and prove lemmas on the functions $\kappa_F$ and $\sigma_F$. 
In Section~2.2 we obtain two lemmas on estimates of $w$. 
In Section~2.3 we obtain sufficient conditions for $\phi\in{\mathcal L}$ to satisfy condition~${\mathcal L}_A$. 
In Section~3 we study the eventual $F$-concavity property of $U_d[\phi]$ with $d\le n/2$ for $\phi\in{\mathcal L}$ 
and prove Theorems~\ref{Theorem:1.2} and \ref{Theorem:1.3}. 
Furthermore, we prove Theorems~\ref{Theorem:1.4} and \ref{Theorem:1.5} using a characterization of log-concavity. 
In Section~4 we study the eventual $F$-concavity property of $U_d[\phi]$ with $d\le n/2$ under condition~(A), and 
prove Theorem~\ref{Theorem:1.6}. 
We also prove Corollaries \ref{Corollary:1.2} and \ref{Corollary:1.3}.
\section{Preliminaries}
%
Throughout this paper, 
for any $x=(x_1,\dots,x_n)$, $y=(y_1,\dots,y_n)\in{\mathbb R}^n$, 
we denote by $\langle x,y\rangle$ 
the standard inner product of $x$ and $y$, that is, 
$\langle x,y\rangle:=\sum_{i=1}^n x_iy_i$. 
Let $e_1$, \dots, $e_n$ be the standard basis of ${\mathbb R}^n$. 
For any measurable set $E$ in ${\mathbb R}^n$, 
we denote by $\chi_E$ and $|E|$ the characteristic function of $E$ and the $n$-dimensional Lebesgue measure of $E$, respectively. 
We also use $C$ to denote generic positive constants  
and may take different values within a calculation. 
%
\subsection{$F$-concavity}
We recall some properties of $F$-concavity. 
Lemmas~\ref{Lemma:2.1} and \ref{Lemma:2.2} concern with the strength of $F$-concavity. 
See \cite{IST05}*{Lemmas~2.4 and 2.5}.
\begin{lemma}
\label{Lemma:2.1}
Let $F_1$ and $F_2$ be admissible on $I=[0,a)$ with $a\in(0,\infty]$. 
Let $f_{F_1}$ {\rm ({\it resp.}~$f_{F_2}$)} be the inverse function of $F_1$ in $J_{F_1}:=F_1({\rm int}\,I)$ {\rm ({\it resp.}~$F_2$ in $J_{F_2}:=F_2({\rm int}\,I)$)}.  
Then $F_2$-concavity is stronger than $F_1$-concavity in ${\mathcal A}(I)$ 
if and only if $F_1(f_{F_2}) $ is concave in $J_{F_2}$  {\rm ({\it or, equivalently, 
$F_2(f_{F_1}) $ is convex in $J_{F_1}$})}.
\end{lemma}
\begin{lemma}
\label{Lemma:2.2}
Let $F_1$ and $F_2$ be admissible on $I=[0,a)$ with $a\in(0,\infty]$. 
Then ${\mathcal C}[F_1]={\mathcal C}[F_2]$ if and only if 
there exists a pair $(A,B)\in (0,\infty)\times{\mathbb R}$ such that 
$$
F_1(\tau)=AF_2(\tau)+B\quad\mbox{for}\quad \tau\in {\rm int}\,I.
$$ 
\end{lemma}
Lemma~\ref{Lemma:2.3} characterizes power concavity by the closedness of $F$-concavity with respect to scalar multiplication 
(see \cite{IST05}*{Lemma~2.8}). See also \cite{IST03}*{Theorem~3.3}.
\begin{lemma}
\label{Lemma:2.3}
Let $F$ be admissible on $I=[0,a)$ with $a\in(0,\infty]$. 
Assume that there exists $\varepsilon_*>0$ such that 
the following holds for every $\varepsilon\in(0,\varepsilon_*]$:
if $f\in{\mathcal C}[F]$ and $(1+\varepsilon)f\in{\mathcal A}(I)$, then 
$\kappa f\in{\mathcal C}[F]$ holds for $\kappa\in[(1+\varepsilon)^{-1},1+\varepsilon]$. 
Then 
$$
{\mathcal C}[F]={\mathcal C}[\Phi_\alpha]\cap{\mathcal A}(I)
$$ 
for some $\alpha\in{\mathbb R}$, 
where $\Phi_\alpha$ is as in {\rm (E1)}.
\end{lemma}
We state two lemmas on necessary conditions for the $F$-concavity property to be preserved by the heat flow. 
Lemma~\ref{Lemma:2.4} follows from \cite{IST03}*{Lemma~4.1}. 
\begin{lemma}
\label{Lemma:2.4}
Let $F$ be admissible in $[0,a)$ with $a\in(0,\infty]$. 
Assume that the $F$-concavity property is preserved by the heat flow in ${\mathbb R}^n$. 
Then, for any $\kappa\in(0,a)$, the function 
$[0,\infty)\ni r\mapsto F(\kappa e^{-r^2})$ is concave. 
\end{lemma}
%
\begin{lemma}
\label{Lemma:2.5}
Let $F$ be admissible in $[0,\infty)$. 
Assume that, for any $\kappa>0$, the function 
$[0,\infty)\ni r\mapsto F(\kappa e^{-r^2})$ is concave. 
Then $F$-concavity is weaker than log-concavity in ${\mathcal A}([0,\infty))$. 
\end{lemma}
{\bf Proof.}
Let $F$ be admissible in $[0,\infty)$. 
Assume that, for any $\kappa>0$, the function 
$[0,\infty)\ni r\mapsto F(\kappa e^{-r^2})$ is concave. 
Then we observe from \cite{IST03}*{Lemma~4.2} that the function $F(e^r)$ is concave for $r\in{\mathbb R}$. 
This together with Lemma~\ref{Lemma:2.1} implies that $F$-concavity is weaker than log-concavity in ${\mathcal A}([0,\infty))$. 
Thus Lemma~\ref{Lemma:2.5} follows. 
$\Box$\vspace{5pt}
\newline
Furthermore, by Lemma~\ref{Lemma:2.1} we have:
\begin{lemma}
\label{Lemma:2.6}
Assume that $F_1$ and $F_2$ satisfy condition~{\rm ($\mbox{F}_a$)} with $a\in(0,\infty)$. 
Then 
$\kappa_{F_1}(r)\ge\kappa_{F_2}(r)$ for $r\in(0,\infty)$ if and only if $F_1$-concavity is stronger than $F_2$-concavity in ${\mathcal A}(I)$.
\end{lemma}
{\bf Proof.}
Thanks to Remark~\ref{Remark:1.1}~(1), 
it suffices to consider the case of $a=1$. 
Since
$$
\kappa_{F_i}(r)
=-r\frac{{\mathcal F}_i''(r)}{{\mathcal F}_i'(r)}=-r\frac{e^{-r}F_i'(e^{-r})+e^{-2r}F''_i(e^{-r})}{-e^{-r}F_i'(e^{-r})}
=r+re^{-r}\frac{F_i''(e^{-r})}{F_i'(e^{-r})}
$$
for $r\in(0,\infty)$, where $i=1,2$, we have 
\begin{equation}
\label{eq:2.1}
\mbox{$\kappa_{F_1}(r)\ge\kappa_{F_2}(r)$ for $r\in(0,\infty)$}
\quad\mbox{if and only if}\quad 
\mbox{$\displaystyle{\frac{F''_1(\tau)}{F_1'(\tau)}\ge \frac{F''_2(\tau)}{F_2'(\tau)}}$ for $\tau\in(0,1)$}.
\end{equation}
 
Let $f_{F_1}$ be the inverse function of $F_1$, and set $J_{F_1}:=F_1((0,1))$. 
Since $F_1'>0$ in $(0,1)$ and $F_1(f_{F_1}(r))=r$ for $r\in J_{F_1}$, 
we see that
\begin{equation}
\label{eq:2.2}
f'_{F_1}(r)=\frac{1}{F_1'(f_{F_1}(r))}>0,\quad
F_1''(f_{F_1}(r))f_{F_1}'(r)^2+F_1'(f_{F_1}(r))f_{F_1}''(r)=0,
\end{equation}
for $r\in J_{F_1}$. 
Since $F_2'>0$ in $(0,1)$, 
by \eqref{eq:2.1} and \eqref{eq:2.2} we obtain 
\begin{align*}
F_2(f_{F_1}(r))'' & =F_2''(f_{F_1}(r))f_{F_1}'(r)^2+F_2'(f_{F_1}(r))f_{F_1}''(r)\\
 & =F_2''(f_{F_1}(r))f_{F_1}'(r)^2-\frac{F_2'(f_{F_1}(r))}{F_1'(f_{F_1}(r))}F_1''(f_{F_1}(r))f_{F_1}'(r)^2\\
 & =f_{F_1}'(r)^2 F_2'(f_{F_1}(r))
 \left[\frac{F_2''(f_{F_1}(r))}{F_2'(f_{F_1}(r))}-\frac{F_1''(f_{F_1}(r))}{F_1'(f_{F_1}(r))}\right]\le 0
\end{align*}
for $r\in J_{F_1}$ if and only if $\kappa_{F_1}(r)\ge\kappa_{F_2}(r)$.
This means that $F_2(f_{F_1})$ is concave in $J_{F_1}$ if and only if $\kappa_{F_1}(r)\ge\kappa_{F_2}(r)$.
Thus Lemma~\ref{Lemma:2.6} follows from Lemma~\ref{Lemma:2.1}.
$\Box$
\vspace{5pt}

We also prove a lemma on the function~$\sigma_F$ defined by \eqref{eq:1.7}.
\begin{lemma}
\label{Lemma:2.7}
Assume condition~{\rm ($\mbox{F}_a$)} with $a=1$. 
Set 
\begin{equation}
\label{eq:2.3}
\nu_F(r):=
\left\{
\begin{array}{ll}
\displaystyle{\frac{r}{\kappa_F(r)}\left(\kappa_F(r)-\frac{1}{2}\right)} & \mbox{if}\quad \kappa_F(r)>0,\vspace{5pt}\\
-\infty & \mbox{if}\quad \kappa_F(r)\le 0,
\end{array}
\right.
\end{equation}
for $r\in(0,\infty)$. Let $\sigma_F$ be as in \eqref{eq:1.7} and $k>0$. 
Then, for any $C_*\in{\mathbb R}$,   
  $$
  \limsup_{r\to\infty}\,(\sigma_F(r)-k\log r)\le C_*
  \quad\mbox{if and only if}\quad
  \limsup_{r\to\infty}\,(\nu_F(r)-2k\log r)\le 2C_*.
  $$
\end{lemma}
{\bf Proof.}
Let $k\ge 0$ and $C_*\in{\mathbb R}$. 
Assume that 
\begin{equation}
\label{eq:2.4}
\displaystyle{\limsup_{r\to\infty}}\,(\sigma_F(r)-k\log r)\le C_*.
\end{equation} 
For any $\epsilon>0$, it follows that
$$
\sigma_F(r)=r\left(\kappa_F(r)-\frac{1}{2}\right)\le k\log r+C_*+\epsilon
$$
for large enough $r$. 
Then
$$
2\kappa_F(r)\le 1+\frac{2k}{r}\log r+\frac{2(C_*+\epsilon)}{r}
$$
for large enough $r$. This implies that 
\begin{align*}
\nu_F(r)=r\left(1-\frac{1}{2\kappa_F(r)}\right) & \le r\left(1-\left(1+\frac{2k}{r}\log r+\frac{2(C_*+\epsilon)}{r}\right)^{-1}\right)\\
 & =2k\log r+2(C_*+\epsilon)+O\left(r^{-1}(\log r)^2\right)
\end{align*}
for large enough $r$ with $\kappa_F(r)>0$. Since $\epsilon$ is arbitrary, we obtain 
\begin{equation}
\label{eq:2.5}
\limsup_{r\to\infty}\left(\nu_F(r)-2k\log r\right)\le 2C_*.
\end{equation}

Next, we assume \eqref{eq:2.5}. 
Then, for any $\epsilon>0$, we have 
$$
\nu_F(r)=r\left(1-\frac{1}{2\kappa_F(r)}\right)\le 2k\log r+2(C_*+\epsilon)
$$
for large enough $r$ with $\kappa_F(r)>0$. 
Then  
$$
\frac{1}{2\kappa_F(r)}\ge 1-\frac{2k}{r}\log r-\frac{2(C_*+\epsilon)}{r}
$$
for large enough $r$ with $\kappa_F(r)>0$. 
This implies that
\begin{align*}
\sigma_F(r)=r\left(\kappa_F(r)-\frac{1}{2}\right)
 & \le r\left(\frac{1}{2}\left(1-\frac{2k}{r}\log r-\frac{2(C_*+\epsilon)}{r}\right)^{-1}-\frac{1}{2}\right)\\
 & =k\log r+C_*+\epsilon+O\left(r^{-1}(\log r)^2\right)
\end{align*}
for large enough $r$ with $\kappa_F(r)>0$. Since $\epsilon$ is arbitrary, we obtain \eqref{eq:2.4}. 
Thus Lemma~\ref{Lemma:2.7} follows.
$\Box$
\subsection{Estimates of $w$}
We obtain two lemmas on estimates of the function $w$ given in \eqref{eq:1.12}. 
Lemma~\ref{Lemma:2.8} (resp.~Lemma~\ref{Lemma:2.9}) concerns with $\phi\in{\mathcal L}$ (resp.~$\phi\in{\mathcal L}_A$). 
\begin{lemma}
\label{Lemma:2.8}
Let $\phi\in{\mathcal L}$. Let $w$ be as in \eqref{eq:1.12}.
Then there exists $C_1>0$ such that
\begin{equation}
\label{eq:2.6}
\frac{1}{2t}\,\underset{y\in P_\phi}{\mbox{{\rm ess inf}}}\,\langle x,y\rangle-C_1\le \log w(x,t)
\le \frac{1}{2t}\,\underset{y\in P_\phi}{\mbox{{\rm ess sup}}}\,\langle x,y\rangle+C_1
\end{equation}
for $(x,t)\in{\mathbb R}^n\times(1,\infty)$.
Furthermore, there exists $C_2>0$ such that
\begin{equation*}
\left|\frac{\partial}{\partial\xi} \log w(x,t)\right|\le C_2t^{-1},\quad \left|\frac{\partial^2}{\partial\xi^2}\log w(x,t)\right|\le C_2t^{-2},
\end{equation*}
for $(x,t)\in{\mathbb R}^n\times(0,\infty)$ and $\xi\in{\mathbb S}^{n-1}$. 
\end{lemma}
{\bf Proof.}
It follows from $\phi\in{\mathcal L}$ that 
\begin{equation}
\label{eq:2.7}
w(x,t)\le \|\phi\|_{L^\infty({\mathbb R}^n)}\int_{P_\phi} e^{\frac{\langle x,y\rangle}{2t}}\,dy\le 
C\exp\left(\frac{1}{2t}\,\underset{y\in P_\phi}{\mbox{{\rm ess sup}}}\,\langle x,y\rangle\right)
\end{equation}
for $x\in{\mathbb R}^n$ and $t>0$. 
On the other hand, since $\phi\ge 0$ and $\phi\not=0$ in ${\mathbb R}^n$, 
we find a measurable set $E\subset P_\phi$ and $\epsilon>0$ such that 
$$
|E|>0,\quad \phi(x)\ge\epsilon>0\quad\mbox{for}\quad x\in E. 
$$
This together with $\phi\in{\mathcal L}$ implies that  
$$
w(x,t)\ge \epsilon\int_E e^{\frac{2\langle x,y\rangle-|y|^2}{4t}}\,dy
\ge \epsilon|E|\exp\left(\frac{1}{4t}\,\underset{y\in P_\phi}{\mbox{{\rm ess inf}}}\,\left(2\langle x,y\rangle-|y|^2\right)\right)
\ge C\exp\left(\frac{1}{2t}\,\underset{y\in P_\phi}{\mbox{{\rm ess inf}}}\,\langle x,y\rangle\right)
$$
for $x\in{\mathbb R}^n$ and $t>1$. This together with \eqref{eq:2.7} implies \eqref{eq:2.6}. 
Furthermore, since $\phi\in {\mathcal L}$, 
we see that
\begin{equation*}
\begin{split}
\left|\left(\frac{\partial}{\partial\xi}\log w\right)(x,t)\right| & \le \frac{|\nabla w(x,t)|}{w(x,t)} 
\le \frac{1}{2tw(x,t)}\int_{{\mathbb R}^n} |y|e^{\frac{2\langle x,y\rangle-|y|^2}{4t}}\phi(y)\,dy\le Ct^{-1},\\
\left|\left(\frac{\partial^2}{\partial\xi^2}\log w\right)(x,t)\right|
 & \le \frac{|\nabla^2 w(x,t)|}{w(x,t)}+\frac{|\nabla w(x,t)|^2}{w(x,t)^2}\\
 & \le \frac{1}{4t^2w(x,t)}\int_{{\mathbb R}^n} |y|^2e^{\frac{2\langle x,y\rangle-|y|^2}{4t}}\phi(y)\,dy+Ct^{-2}\le Ct^{-2},
\end{split}
\end{equation*}
for $(x,t)\in{\mathbb R}^n\times(0,\infty)$ and $\xi\in{\mathbb S}^{n-1}$. 
Thus Lemma~\ref{Lemma:2.8} follows.
$\Box$
\begin{lemma}
\label{Lemma:2.9}
Let $\phi\in {\mathcal L}_A$. 
Assume that 
\begin{equation}
\label{eq:2.8}
\,\underset{y\in P_\phi}{\mbox{{\rm ess sup}}}\,\langle y,e_1\rangle=0.
\end{equation}
Let $w$ be as in \eqref{eq:1.12}. 
Then there exists $C_1>0$ such that 
\begin{equation}
\label{eq:2.9}
\log w(x,t)\le -\log\frac{X}{t}+C_1
\end{equation}
for $x=(X,0)$ with $X>0$ and $t>1$. 
Furthermore, there exists $C_2>0$ such that 
\begin{equation}
\label{eq:2.10}
\begin{split}
 & \left|\frac{\partial}{\partial x_1}(\log w)(x,t)\right|\le C_2X^{-1},\quad
0\le\frac{\partial^2}{\partial x_1^2}(\log w)(x,t)\le C_2X^{-2},\\
 & \left|\frac{\partial^2}{\partial x_1\partial x_i}(\log w)(x,t)\right|\le C_2t^{-1}X^{-1},
\end{split}
\end{equation}
for $x=(X,0)$ with $X\ge 2t$ and $t>1$, where $i=2,\dots,n$.
\end{lemma}
{\bf Proof.}
It follows from $\phi\in{\mathcal L}_A$, \eqref{eq:1.12}, and \eqref{eq:2.8} that
$$
\log w(x,t)\le \log\left[C\|\phi\|_{L^\infty({\mathbb R}^n)}\int_{-\infty}^0 e^{\frac{Xy_1}{2t}}\,dy_1\right]\le\log\frac{Ct}{X}
\le -\log\frac{X}{t}+C
$$
for $x=(X,0)$ with $X>0$ and $t\in(1,\infty)$. This implies \eqref{eq:2.9}.

We apply Jensen's inequality and condition~(A) with \eqref{eq:2.8} to obtain 
\begin{align*}
\left|\frac{\partial}{\partial x_1}(\log w)(x,t)\right| & \le\frac{1}{2t w(x,t)}\int_{{\mathbb R}^n}|y_1|e^{\frac{2Xy_1-|y|^2}{4t}}\phi(y)\,dy\\
 & \le\frac{1}{2t}\left(\frac{1}{w(x,t)}\int_{{\mathbb R}^n}|y_1|^2e^{\frac{2Xy_1-|y|^2}{4t}}\phi(y)\,dy\right)^{1/2}\\
 & \le\frac{C}{t}\left(\int_{{\mathbb R}^n}|y_1|^2e^{\frac{Xy_1}{2t}}\phi(y)\,dy\biggr/\int_{{\mathbb R}^n}e^{\frac{Xy_1}{2t}}\phi(y)\,dy\right)^{1/2}
 \le CX^{-1},
\end{align*}
\begin{equation*}
\begin{split}
 & \frac{\partial^2}{\partial x_1^2}(\log w)(x,t)\\
 & =\frac{1}{4t^2}
\left[\frac{1}{w(x,t)}\int_{{\mathbb R}^n} |y_1|^2e^{\frac{2Xy_1-|y|^2}{4t}}\phi(y)\,dy
-\left(\frac{1}{w(x,t)}\int_{{\mathbb R}^n}y_1e^{\frac{2Xy_1-|y|^2}{4t}}\phi(y)\,dy\right)^2\right]\ge 0,
\end{split}
\end{equation*}
for $x=(X,0)$ with $X\ge 2t$ and $t>1$.
Furthermore, 
$$
0\le\frac{\partial^2}{\partial x_1^2}(\log w)(x,t)\le\frac{1}{4t^2 w(x,t)}\int_{{\mathbb R}^n}|y_1|^2e^{\frac{2Xy_1-|y|^2}{4t}}\phi(y)\,dy
 \le CX^{-2}
$$
for $x=(X,0)$ with $X\ge 2t$ and $t>1$. 
Similarly, we have
\begin{align*}
 & \left|\frac{\partial^2}{\partial x_1\partial x_i}(\log w)(x,t)\right|\\
 & \le \frac{1}{4t^2w(x,t)}\int_{{\mathbb R}^n}|y_1||y_i|e^{\frac{2Xy_1-|y|^2}{4t}}\phi(y)\,dy\\
 & \qquad\quad
 +\frac{1}{4t^2w(x,t)^2}\left(\int_{{\mathbb R}^n}|y_1|e^{\frac{2Xy_1-|y|^2}{4t}}\phi(y)\,dy\right)\left(\int_{{\mathbb R}^n}|y_i|e^{\frac{2Xy_1-|y|^2}{4t}}\phi(y)\,dy\right)\\
 & \le \frac{C}{t^2w(x,t)}\int_{{\mathbb R}^n}|y_1|e^{\frac{2Xy_1-|y|^2}{4t}}\phi(y)\,dy\\
 & \le \frac{C}{t^2}\left(\frac{1}{w(x,t)}\int_{{\mathbb R}^n}|y_1|^2e^{\frac{2Xy_1-|y|^2}{4t}}\phi(y)\,dy\right)^{1/2}
 \le Ct^{-1}X^{-1}
\end{align*}
for $x=(X,0)$ with $X\ge 2t$ and $t>1$.
Thus \eqref{eq:2.10} holds, and Lemma~\ref{Lemma:2.9} follows.
$\Box$
\subsection{Subclass~${\mathcal L}_A$}
We give sufficient conditions for $\phi\in{\mathcal L}$ to satisfy condition~(A). 
\begin{lemma}
\label{Lemma:2.10}
Let $\phi\in{\mathcal L}$ satisfy the following condition.
\begin{itemize}
\item[{\rm (A1)}]
 There exist $\alpha$, $\delta>0$ such that
  $P_\phi\cap B(x_*,\delta)$ is open and $\phi\chi_{B(x_*,\delta)}$ is $\alpha$-concave in ${\mathbb R}^n$
   for all $x_*\in \partial_* P_\phi$, where 
  $$
  \partial_*P_\phi:=\biggr\{x\in\partial P_\phi\,:\,
  \sup_{y\in P_\phi}\langle y-x,\xi\rangle=0\mbox{ for some $\xi\in{\mathbb S}^{n-1}$}\biggr\}.
  $$
\end{itemize}
Then $\phi\in {\mathcal L}_A$. 
\end{lemma}
{\bf Proof.}
Let $\phi\in{\mathcal L}$ satisfy condition~(A1). 
Let $(\xi,z)\in{\mathbb S}^{n-1}\times{\mathbb R}^n$ satisfy $\sup_{y\in P_\phi}\langle y-z,\xi\rangle=0$.  
We can assume, without loss of generality, that $\xi=e_1$, $z=0$, 
\begin{equation}
\label{eq:2.11}
P_\phi\subset\{(x_1,x')\,:\,x_1\le 0,\,x'\in{\mathbb R}^{n-1}\},\quad\mbox{and}\quad
\partial P_\phi\cap \Pi\not=\emptyset,
\end{equation}
where $\Pi:=\{0\}\times{\mathbb R}^{n-1}$. 
For the proof, it suffices to prove that 
\begin{equation}
\label{eq:2.12}
 \int_{{\mathbb R}^n}y_1^2e^{ky_1}\phi(y)\,dy
 \le Ck^{-2}\int_{{\mathbb R}^n} e^{ky_1}\phi(y)\,dy\quad\mbox{for}\quad k\ge 1.
\end{equation}

Let $\alpha$, $\delta>0$ be as in condition~(A1). 
Thanks to the compactness of $\partial P_\phi$, 
by \eqref{eq:2.11} 
we find finite points $\{x_j\}_{j=1}^m\subset\partial P_\phi\cap\Pi$
such that 
$$
\partial P_\phi\cap\Pi\subset\bigcup_{j=1}^m B(x_j,\delta).
$$
Furthermore, 
there exists $\ell>0$ such that 
\begin{equation}
\label{eq:2.13}
P_\phi^\ell:=\{(y_1,y)\in P_\phi\,:\,y_1>-\ell\}\subset \bigcup_{j=1}^m B(x_j,\delta).
\end{equation}
On the other hand, 
it follows from condition~(A1) that, for any $j=1,\dots,m$, 
$\phi\chi_{B(x_j,\delta)}$ is $\alpha$-concave in ${\mathbb R}^n$. 
Then we observe that 
\begin{equation}
\label{eq:2.14}
\mbox{$\phi^\alpha$ is concave in $P_\phi\cap B(x_j,\delta)$ and $P_\phi\cap B(x_j,\delta)$ is convex and open}. 
\end{equation}
In particular, $\phi$ is continuous in $P_\phi\cap B(x_j,\delta)$. 
Then we find an (non-empty) open set $E$ in ${\mathbb R}^n$ such that 
$$
E\subset P_\phi^{\ell/2},\qquad \underset{x\in E}{\mbox{inf}}\,\phi(x)>0. 
$$
Then we have
$$
\int_{P_\phi} e^{ky_1}\phi(y)\,dy\ge C\,\left(\underset{x\in E}{\mbox{ess inf}}\,\phi(x)\right)\,\int_E e^{ky_1}\,dy
\ge Ce^{-\frac{1}{2}k\ell}.
$$
This implies that
\begin{equation}
\label{eq:2.15}
\int_{P_\phi\setminus P_\phi^\ell} e^{\frac{1}{2} ky_1}\phi(y)\,dy
\le C\|\phi\|_{L^\infty({\mathbb R}^n)}\int_{P_\phi\setminus P_\phi^\ell} e^{-\frac{1}{2}k\ell}\,dy
\le Ce^{-\frac{1}{2}k\ell}\le C\int_{P_\phi} e^{ky_1}\phi(y)\,dy.
\end{equation}
On the other hand, 
for any fixed $j=1,\dots,m$, 
we observe from $x_j\in\Pi$ that
$$
z:=\frac{x_j+y}{2}\in P_\phi\cap B(x_i,\delta),\quad
\frac{1}{2}\phi(y)^\alpha=\frac{1}{2}\left[\phi(y)^\alpha+\phi(x_j)^\alpha\right]
\le\phi^\alpha\left(\frac{x_j+y}{2}\right)=\phi^\alpha(z),
$$
for $y\in P_\phi\cap B(x_j,\delta)$. 
These imply that 
\begin{equation}
\label{eq:2.17}
\int_{P_\phi\cap B(x_j,\delta)} e^{\frac{1}{2}ky_1}\phi(y)\,dy
\le C\int_{P_\phi\cap B(x_j,\delta)} e^{kz_1}\phi(z)\,dz,\quad j=1,\dots, m.
\end{equation}
Combining \eqref{eq:2.13}, \eqref{eq:2.15}, and \eqref{eq:2.17}, 
we see that 
\begin{equation*}
\begin{split}
 \int_{{\mathbb R}^n}y_1^2e^{ky_1}\phi(y)\,dy & = k^{-2}\int_{P_\phi}(ky_1)^2e^{ky_1}\phi(y)\,dy
 \le Ck^{-2}\int_{P_\phi} e^{\frac{1}{2}ky_1}\phi(y)\,dy\\
 & =Ck^{-2}\left(\int_{P_\phi\setminus P_\phi^\ell}+\sum_{j=1}^m \int_{P_\phi\cap B(x_j,\delta)}\right) e^{\frac{1}{2}ky_1}\phi(y)\,dy\\
 & \le Ck^{-2}\int_{P_\phi} e^{ky_1}\phi(y)\,dy+Ck^{-2}\sum_{j=1}^m \int_{P_\phi\cap B(x_j,\delta)} e^{kz_1}\phi(z)\,dz\\
 & \le Ck^{-2}\int_{{\mathbb R}^n} e^{ky_1}\phi(y)\,dy. 
\end{split}
\end{equation*}
This implies \eqref{eq:2.12}. Thus Lemma~\ref{Lemma:2.10} follows.
$\Box$
\vspace{3pt}
\newline
As an application of Lemma~\ref{Lemma:2.10}, we see that 
$\chi_\Omega\in{\mathcal L}_A$ if $\Omega$ is a bounded smooth open set in~${\mathbb R}^n$.  
Furthermore, 
we easily obtain the following properties.
\begin{itemize}
  \item[(A2)] 
  Let $\phi\in{\mathcal L}_A$ and $\tilde{\phi}\in{\mathcal L}$. 
  If there exists $C\ge 1$ such that $C^{-1}\phi(x)\le\tilde{\phi}(x)\le C\phi(x)$ for almost all $x\in{\mathbb R}^n$, 
  then $\tilde{\phi}\in{\mathcal L}_A$. 
  \item[(A3)] 
  If $\phi$, $\tilde{\phi}\in{\mathcal L}_A$, then $c_1\phi+c_2\tilde{\phi}\in{\mathcal L}_A$ for $c_1$, $c_2>0$.  
\end{itemize}
\section{Eventual $F$-concavity in ${\mathcal L}$}
In this section we study the eventual $F$-concavity property of $U_d[\phi]$ for $\phi\in{\mathcal L}$, 
and prove Theorems~\ref{Theorem:1.2}--\ref{Theorem:1.5} and Corollary~\ref{Corollary:1.1}.

We first prove Theorems~\ref{Theorem:1.2} and \ref{Theorem:1.3}. 
Let $F$ satisfy condition~($\mbox{F}_a$) with $a\in(0,\infty)$, $d\le n/2$, and $\phi\in{\mathcal L}$. Assume that $M_\phi<a$ if $d=n/2$. 
Thanks to Remark~\ref{Remark:1.1}~(1), we can assume, without loss of generality, that $a=1$. 
\begin{lemma}
\label{Lemma:3.1}
Assume condition~{\rm ($\mbox{F}_a$)} with $a=1$. Let $\phi\in{\mathcal L}$. 
Let $v_d$ be as in \eqref{eq:1.11}. 
\begin{itemize}
  \item[{\rm (1)}] 
  Let $d<n/2$. Then 
  \begin{equation}
  \label{eq:3.1}
  \lim_{t\to\infty}\inf_{x\in{\mathbb R}^n}v_d(x,t)=\infty.
  \end{equation}
  \item[{\rm (2)}] 
  Let $d=n/2$ and $M_\phi<1$. 
  Then 
 \begin{equation}
  \label{eq:3.2}
  v_d(x,t)\ge-\log M_\phi>0
 \end{equation}
  for $(x,t)\in{\mathbb R}^n\times(0,\infty)$. 
  Furthermore, 
  for any $m>0$, 
  there exists $L>0$ such that 
 \begin{equation}
  \label{eq:3.3}
  v_d(x,t)\ge m
 \end{equation}
  for $x\in{\mathbb R}^n\setminus B(0,L\sqrt{t})$ and $t\in[1,\infty)$. 
\end{itemize}
\end{lemma}
{\bf Proof.}
It follows from \eqref{eq:1.1} that 
$$
U_d(x,t)\le (4\pi t)^{d-\frac{n}{2}}M_\phi
$$ 
for $(x,t)\in{\mathbb R}^n\times(0,\infty)$. 
This implies that 
$$
\inf_{x\in{\mathbb R}^n} v_d(x,t)=-\log\left(\sup_{x\in{\mathbb R}^n} U_d(x,t)\right)\ge \left(\frac{n}{2}-d\right)\log(4\pi t)-\log M_\phi
$$
for $(x,t)\in{\mathbb R}^n\times(0,\infty)$. This implies \eqref{eq:3.1} and \eqref{eq:3.2}. 
Furthermore, in the case of $d=n/2$, 
by Lemma~\ref{Lemma:2.8} we have 
$$
v_d(x,t)=\frac{|x|^2}{4t}-\log w\ge \frac{|x|^2}{4t}-C\frac{|x|}{t}-C
=\frac{|x|}{t}\left(\frac{|x|}{4}-C\right)-C
$$
for $(x,t)\in{\mathbb R}^n\times[1,\infty)$.
Then, taking large enough $L>0$, we obtain 
$$
v_d(x,t)
\ge \frac{L^2}{8}-C\ge\frac{L^2}{16}
$$
for $x\in{\mathbb R}^n\setminus B(0,L\sqrt{t})$ and $t\in[1,\infty)$. 
This implies \eqref{eq:3.3}, and Lemma~\ref{Lemma:3.1} follows.
$\Box$
\vspace{5pt}

Assume condition~($\mbox{F}_a$) with $a=1$. 
Since ${\mathcal F}(r)$ is defined for $r>0$, 
by Lemma~\ref{Lemma:3.1}
we find $T>0$ such that the function 
$$
F\left(U_d[\phi](x,t)\right)={\mathcal F}\left(-\log\left(U_d[\phi](x,t)\right)\right)
={\mathcal F}(v_d(x,t))
$$
is well-defined for $(x,t)\in{\mathbb R}^n\times[T,\infty)$. 
For any $\xi\in{\mathbb S}^{n-1}$, it follows from \eqref{eq:1.3} and \eqref{eq:1.11} that
\begin{equation*}
\begin{split}
\frac{\partial^2}{\partial\xi^2}F\left(U_d[\phi](x,t)\right)
 & ={\mathcal F}''(v_d(x,t))\left(\frac{\partial}{\partial\xi} v_d(x,t)\right)^2+{\mathcal F}'(v_d(x,t))\frac{\partial^2}{\partial\xi^2}v_d(x,t)\\
 & =\frac{{\mathcal F}'(v_d(x,t))}{v_d(x,t)}I_F[v_d](x,t;\xi),
\end{split}
\end{equation*}
where 
\begin{equation}
\label{eq:3.4}
\begin{split}
 & I_F[v_d](x,t;\xi)\\
 & :=-\kappa_F(v_d(x,t))\left(\frac{\partial}{\partial\xi} v_d(x,t)\right)^2+v_d(x,t)\frac{\partial^2}{\partial\xi^2}v_d(x,t)\\
 & =-\kappa_F(v_d(x,t))\left(\frac{\langle x,\xi\rangle}{2t}-\left(\frac{\partial}{\partial\xi}\log w\right)(x,t)\right)^2\\
 & \qquad
 +\left(\left(\frac{n}{2}-d\right)\log(4\pi t)+\frac{|x|^2}{4t}-\log w(x,t)\right)\left(\frac{1}{2t}-\left(\frac{\partial^2}{\partial\xi^2}\log w\right)(x,t)\right).
\end{split}
\end{equation}
Since ${\mathcal F}'<0$ in $(0,\infty)$, we see that 
\begin{equation}
\label{eq:3.5}
\frac{\partial^2}{\partial\xi^2}F\left(U_d[\phi](x,t)\right)\le 0\quad\mbox{if and only if}\quad I_F[v_d](x,t;\xi)\ge 0
\end{equation}
for $(x,t)\in{\mathbb R}^n\times[T,\infty)$.

\begin{proposition}
\label{Proposition:3.1}
Assume condition~{\rm ($\mbox{F}_a$)} with $a=1$. 
Let $\phi\in{\mathcal L}$ and $d\le n/2$. 
Assume that $M_\phi<1$ if $d=n/2$. 
\begin{itemize}
  \item[{\rm (1)}] 
  Let $L>0$. There exists $T_1>0$ such that 
  \begin{equation}
  \label{eq:3.6}
  \frac{\partial^2}{\partial\xi^2}F\left(U_d[\phi](x,t)\right)\le 0,
  \quad \xi\in{\mathbb S}^{n-1},
  \end{equation}
  for $x\in B(0,L\sqrt{t})$ and $t\in[T_1,\infty)$ with $\kappa_F(v_d(x,t))\le 1/2$.
  \item[{\rm (2)}] 
  Let $\kappa_*<1/2$. There exists $T_2>0$ such that \eqref{eq:3.6} holds
   for $(x,t)\in {\mathbb R}^n\times[T_2,\infty)$ with $\kappa_F(v_d(x,t))\le\kappa_*$ and $\xi\in{\mathbb S}^{n-1}$.
\end{itemize}
\end{proposition}
{\bf Proof.}
Let $\phi\in{\mathcal L}$. 
We prove assertion~(1). Let $L>0$ and $\xi\in{\mathbb S}^{n-1}$. 
We can assume, without loss of generality, that 
\begin{equation}
\label{eq:3.7}
\int_{{\mathbb R}^n}y_i\phi(y)\,dy=0,\quad i=1,\dots,n. 
\end{equation}
Then 
\begin{equation}
\label{eq:3.8}
\lim_{t\to\infty}w(x,t)=\lim_{t\to\infty}\int_{{\mathbb R}^n}e^{\frac{2\langle x,y\rangle-|y|^2}{4t}}\,\phi(y)\,dy
=M_\phi.
\end{equation}
Furthermore, 
\begin{equation}
\label{eq:3.9}
\begin{split}
 & \lim_{t\to\infty}t\left(\frac{\partial}{\partial\xi}w\right)(x,t)
=\frac{1}{2}\lim_{t\to\infty}\int_{{\mathbb R}^n}\langle \xi,y\rangle e^{\frac{2\langle x,y\rangle-|y|^2}{4t}}\,\phi(y)\,dy
=\int_{{\mathbb R}^n}\langle \xi,y\rangle\,\phi(y)\,dy=0,\\
 & \lim_{t\to\infty}t^2\left(\frac{\partial^2}{\partial\xi^2}w\right)(x,t)
=\frac{1}{4}\lim_{t\to\infty}\int_{{\mathbb R}^n}\langle \xi,y\rangle^2 e^{\frac{2\langle x,y\rangle-|y|^2}{4t}}\,\phi(y)\,dy
=\int_{{\mathbb R}^n}\langle \xi,y\rangle^2\,\phi(y)\,dy>0,
\end{split}
\end{equation}
uniformly for $x\in B(0,L\sqrt{t})$. 
Since $M_\phi<1$ if $d=n/2$, 
by \eqref{eq:3.4}, \eqref{eq:3.8}, and \eqref{eq:3.9} we have 
\begin{equation*}
\begin{split}
 & I_F[v_d](x,t;\xi)\\
 & \ge -\frac{1}{2}\left(\frac{\langle x,\xi\rangle}{2t}+o(t^{-1})\right)^2\\
 & \qquad\quad
 +\left(\left(\frac{n}{2}-d\right)\log(4\pi t)+\frac{|x|^2}{4t}-(\log M_\phi)(1+o(1))\right)\left(\frac{1}{2t}+O(t^{-2})\right)\\
 & =\frac{1}{2}\left(\frac{n}{2}-d\right)t^{-1}\log(4\pi t)
 +\frac{1}{8t^2}(|x|^2-\langle x,\xi\rangle^2)-\frac{\log M_\phi}{2t}(1+o(1))+o(t^{-\frac{3}{2}})\\
 & \ge\frac{1}{2}\left(\frac{n}{2}-d\right)t^{-1}\log(4\pi t)-\frac{\log M_\phi}{2t}(1+o(1))+o(t^{-\frac{3}{2}})>0
\end{split}
\end{equation*}
for $x\in B(0,L\sqrt{t})$ and large enough $t$ with $\kappa_F(v_d(x,t))\le 1/2$. 
This together with \eqref{eq:3.5} implies assertion~(1). 

We prove assertion~(2). 
Let $\kappa_*<1/2$. Let $L>0$ be large enough. 
By assertion~(1) it suffices to consider the case where $x\in{\mathbb R}^n\setminus B(0,L\sqrt{t})$ and large enough $t$ with $\kappa_F(v_d(x,t))\le\kappa_*$. 
By Lemma~\ref{Lemma:2.8} and \eqref{eq:3.4} we have
\begin{equation}
\label{eq:3.10}
\begin{split}
I_F[v_d](x,t;\xi) & \ge
 -\kappa_*\left(\frac{\langle x,\xi\rangle}{2t}+O(t^{-1})\right)^2\\
 & \qquad
 +\left(\left(\frac{n}{2}-d\right)\log(4\pi t)+\frac{|x|^2}{4t}+O(t^{-1}|x|)-C\right)\left(\frac{1}{2t}+O(t^{-2})\right)\\
 & =\frac{1}{2}\left(\frac{n}{2}-d\right)t^{-1}(\log 4\pi t)-\frac{\kappa_*}{4t^2}\langle x,\xi\rangle^2+\frac{|x|^2}{8t^2}-\frac{C}{2t}\\
 & \qquad
 +O(t^{-3}|x|^2)+O(t^{-2}|x|)+O(t^{-2}\log t)\\
 & \ge \left(\frac{1}{2}-\kappa_*\right)\frac{|x|^2}{4t^2}-\frac{C}{2t}+O(t^{-3}|x|^2)+O(t^{-2}|x|)+O(t^{-2}\log t)\\
 & \ge\frac{1}{2}\left(\frac{1}{2}-\kappa_*\right)\frac{|x|^2}{4t^2}(1+o(1))+\frac{1}{2}\left(\frac{1}{2}-\kappa_*\right)\frac{L^2}{4t}-\frac{C}{2t}>0
\end{split}
\end{equation}
for $x\in {\mathbb R}^n\setminus B(0,L\sqrt{t})$ and large enough $t$ with $\kappa_F(v_d(x,t))\le \kappa_*$. 
This together with assertion~(1) and \eqref{eq:3.5} implies assertion~(2). 
Thus Proposition~\ref{Proposition:3.1} follows.
$\Box$\vspace{5pt}

Now we are ready to prove Theorems~\ref{Theorem:1.2} and \ref{Theorem:1.3}. 
It suffices to consider the case of $a=1$. 
\vspace{5pt}
\newline
{\bf Proof of  Theorem~\ref{Theorem:1.2}.}
Let $d<n/2$, $\phi\in{\mathcal L}$, and $a=1$. 
We first prove assertion~(1). 
Assume that $\kappa_F^*<1/2$. 
Then we find $R>0$ and $\kappa'>0$ such that 
\begin{equation}
\label{eq:3.11}
\kappa_F(r)\le\kappa'<\frac{1}{2},\quad r\in(R,\infty). 
\end{equation}
This together with Lemma~\ref{Lemma:3.1}~(1) implies that 
$$
\kappa_F(v_d(x,t))\le\kappa'<\frac{1}{2}
$$
for $x\in{\mathbb R}^n$ and large enough $t$. 
Then, by Proposition~\ref{Proposition:3.1}~(2) we find $T>0$ such that  
$$
\frac{\partial^2}{\partial\xi^2}F\left(U_d[\phi](x,t)\right)\le 0,
\quad \xi\in{\mathbb S}^{n-1},
$$
for $(x,t)\in{\mathbb R}^n\times[T,\infty)$. 
This implies that $U_d[\phi]$ is $F$-concave in ${\mathbb R}^n$ for $t\in[T,\infty)$. 
Thus assertion~(1) follows. 

Next, we prove assertion~(2). 
Assume that $\kappa_F^*>1/2$. 
Similarly to \eqref{eq:3.10}, 
by Lemmas~\ref{Lemma:2.8} and \ref{Lemma:3.1} we have  
\begin{align*}
I_F[v_d](te_1,t;e_1) & =-\kappa(v_d(te_1,t))\left(\frac{1}{2}+O(t^{-1})\right)^2\\
 & \qquad\quad
 +\left(\left(\frac{n}{2}-d\right)\log(4\pi t)+\frac{t}{4}+O(1)\right)\left(\frac{1}{2t}+O(t^{-2})\right)\\
 & =\frac{1}{2}\left(\frac{n}{2}-d\right)t^{-1}(\log 4\pi t)(1+o(1))+\frac{1}{4}\left(\frac{1}{2}-\kappa(v_d(te_1,t))(1+O(t^{-1}))\right)
\end{align*}
for large enough $t$. 
Since $\kappa_F^*>1/2$,  we find a sequence $\{r_j\}\subset{\mathbb R}$ and $\kappa''\in{\mathbb R}$ such that 
$$
\lim_{j\to\infty}r_j=\infty,\qquad \kappa(r_j)>\kappa''>\frac{1}{2}\quad\mbox{for $j=1,2,\dots$}.
$$
By Lemma~\ref{Lemma:3.1}, applying by the intermediate value theorem, 
for any large enough $T>0$, 
we find $t_T\in[T,\infty)$ such that $v_d(t_Te_1,t_T)=r_j$ for some $j=1,2,\dots$. 
Then, taking large enough $T$ if necessary, we see that 
$$
I_F[v_d](t_Te_1,t_T;e_1)\le\left(\frac{n}{2}-d\right)t_T^{-1}(\log 4\pi t_T)+\frac{1}{4}\left(\frac{1}{2}-\kappa''(1+O(T^{-1}))\right)<0. 
$$
This together with \eqref{eq:3.5}  implies 
that $U_d[\phi](\cdot,t_T)$ is not $F$-concave in ${\mathbb R}^n$. 
Thus assertion~(2) follows, and the proof of Theorem~\ref{Theorem:1.2} is complete. 
$\Box$
\vspace{5pt}

\noindent
{\bf Proof of  Theorem~\ref{Theorem:1.3}.}
Let $d=n/2$, $\phi\in{\mathcal L}$ wtih $M_\phi<1$, and $a=1$. 
We prove assertion~(1).  Assume \eqref{eq:1.6}. 
Let $L>0$ be large enough. 
Then, by Proposition~\ref{Proposition:3.1}~(1)
we find $T>0$ such that 
\begin{equation}
\label{eq:3.12}
\frac{\partial^2}{\partial\xi^2}F\left(U_d[\phi](x,t)\right)\le 0,
\quad \xi\in{\mathbb S}^{n-1},
\end{equation}
for $x\in B(0,L\sqrt{t})$ and $t\in[T,\infty)$. 

By \eqref{eq:1.6} we find $R>0$ and $\kappa'>0$ satisfying \eqref{eq:3.11}. 
By Lemma~\ref{Lemma:3.1}~(2), 
taking large enough $L$ if necessary, we have
$$
v_d(x,t)\ge R,\quad 
\kappa_F[v_d(x,t)]\le\kappa'<\frac{1}{2},
$$
for $x\in{\mathbb R}^n\setminus B(0,L\sqrt{t})$ and $t\in[1,\infty)$. 
Then, applying Proposition~\ref{Proposition:3.1}~(2) and taking large enough~$T'\in[T,\infty)$, 
we deduce that 
\begin{equation}
\label{eq:3.13}
\frac{\partial^2}{\partial\xi^2}F\left(U_d[\phi](x,t)\right)\le 0,
\quad \xi\in{\mathbb S}^{n-1},
\end{equation}
for $x\in {\mathbb R}^n\setminus B(0,L\sqrt{t})$ and $t\in[T',\infty)$. 
Combining \eqref{eq:3.12} and \eqref{eq:3.13}, we observe that 
$$
\frac{\partial^2}{\partial\xi^2}F\left(U_d[\phi](x,t)\right)\le 0,
\quad \xi\in{\mathbb S}^{n-1},
$$
for $x\in {\mathbb R}^n$ and $t\in[T',\infty)$. 
Thus $U_d[\phi]$ possesses the eventual $F$-concavity property, and assertion~(1) follows. 

Next, we prove assertion~(2). 
Assume that $\kappa_F(r_*)>1/2$ for some $r_*\in(0,\infty)$. 
Let $\phi\in{\mathcal L}$ be such that $M_\phi\in(e^{-r_*/2},1)$. 
We can assume, without loss of generality, that $\phi$ satisfies \eqref{eq:3.7}. 
For any $L>0$, 
by \eqref{eq:1.11} and \eqref{eq:3.8} we have
$$
\lim_{t\to\infty}v_d(r\sqrt{t}e_1,t)=\frac{r^2}{4}-\log M_\phi
$$
uniformly for $r\in(0,L)$. 
Then, for any large enough $t>0$, 
applying the intermediate value theorem, we find $r(t)>0$ such that 
\begin{equation}
\label{eq:3.14}
v_d(r(t)\sqrt{t}e_1,t)=r_*,\quad \lim_{t\to\infty}r(t)=2\sqrt{\log M_\phi+r_*}>2\sqrt{r_*/2}>0. 
\end{equation}
It follows from \eqref{eq:3.4}, \eqref{eq:3.8}, \eqref{eq:3.9}, and \eqref{eq:3.14} that 
\begin{equation*}
\begin{split}
 & I_F[v_{n/2}](r(t)\sqrt{t}e_1,t;e_1)\\
 & =-\kappa_F(r_*)\left(\frac{r(t)}{2\sqrt{t}}+o(t^{-1})\right)^2
 +\left(\frac{r(t)^2}{4}-(\log M_\phi)(1+o(1))\right)\left(\frac{1}{2t}+O(t^{-2})\right)\\
 & =\frac{1+o(1)}{2t}|\log M_\phi|+\left(\frac{1}{2}-\kappa_F(r_*)\right)\frac{r(t)^2}{4t}+o(t^{-\frac{3}{2}})\\
 & \le\frac{1}{t}|\log M_\phi|+\left(\frac{1}{2}-\kappa_F(r_*)\right)\frac{r_*}{2t}
\end{split}
\end{equation*}
for large enough $t$. Taking $M_\phi$ close enough to $1$ if necessary, we deduce that 
$$
I_F[v_{n/2}](r(t)\sqrt{t}e_1,t;e_1)<0
$$
for large enough $t$. 
This implies that $U_{n/2}[\phi](t)$ is not $F$-concave in ${\mathbb R}^n$ for large enough $t$. 
Thus assertion~(2) follows, and the proof of Theorem~\ref{Theorem:1.3} is complete. 
$\Box$
\vspace{5pt}
\newline
By Theorems~\ref{Theorem:1.2} and \ref{Theorem:1.3} 
we prove Corollary~\ref{Corollary:1.1}.  
\vspace{5pt}
\newline
{\bf Proof of Corollary~\ref{Corollary:1.1}.}
Let us consider $\beta$-log-concavity, where $\beta>0$, that is, $F=\Psi_\beta$ (see (E2)). 
It follows from Remark~\ref{Remark:1.1}~(2) that $\kappa_F(r)=-\beta+1$ for $r\in(0,\infty)$. 
Then Corollary~\ref{Corollary:1.1} follows from Theorems~\ref{Theorem:1.2} and \ref{Theorem:1.3}. 
$\Box$
\vspace{5pt}

Next, we develop the arguments in \cite{IST03}*{Section~4} to prove Theorem~\ref{Theorem:1.5}. 
\vspace{5pt}
\newline
{\bf Proof of Theorem~\ref{Theorem:1.5}.}
Let $d>n/2$ and $\phi\in{\mathcal L}$.
 Assume that $U_d[\phi]$ possesses the eventual $F$-concavity property.
Let $k>0$, and set 
$$
m:=k/M_\phi,
\quad
\eta(t):=\left(d-\frac{n}{2}\right)\log(4\pi t)-\log m.
$$
For any $r\in[0,\infty)$, set $x(t:r):=2\sqrt{t}\sqrt{\eta(t)+r^2}\,e_1$
for large enough $t$.
Since 
$$
\exp\left(-\frac{|x(t:r)|^2}{4t}\right)=(4\pi t)^{-d+\frac{n}{2}}me^{-r^2},
$$
we have 
\begin{equation}
\label{eq:3.15}
\lim_{t\to\infty}U_d(x(t:r),t)=me^{-r^2}\lim_{t\to\infty}\int_{{\mathbb R}^n}e^{\frac{\langle x(t:r),y\rangle-|y|^2}{4t}}\phi(y)\,dy
=me^{-r^2}\int_{{\mathbb R}^n}\phi(y)\,dy=ke^{-r^2}
\end{equation}
for $r\in[0,\infty)$.
Let $r_0$, $r_1\in[0,\infty)$ and $\mu\in(0,1)$. 
Set 
$$
x_0(t):=x(t,r_0),\quad x_1(t):=x(t,r_1),\quad x_\mu(t):=(1-\mu)x_0(t)+\mu x_1(t),
$$
for large enough $t$. 
Since $U_d$ possesses the eventual $F$-concavity property, we have 
\begin{equation}
\label{eq:3.16}
F(U_d(x_\mu(t),t))\ge(1-\mu)F(U_d(x_0(t),t))+\mu F(U_d(x_1(t),t))
\end{equation}
for large enough $t$. It follows from \eqref{eq:3.15} that 
\begin{equation}
\label{eq:3.17}
\lim_{t\to\infty}\left[(1-\mu)F(U_d(x_0(t),t))+\mu F(U_d(x_1(t),t))\right]
=(1-\mu)F\left(ke^{-r_0^2}\right)+\mu F\left(e^{-r_1^2}\right).
\end{equation}
On the other hand, it follows that 
\begin{align*}
\frac{|x_\mu(t)|}{2\sqrt{t}}
 & =\sqrt{\eta(t)}\left((1-\mu)\sqrt{1+\frac{r_0^2}{\eta(t)}}+\mu\sqrt{1+\frac{r_1^2}{\eta(t)}}\right)\\
 & =\sqrt{\eta(t)}\left(1+\frac{(1-\mu)r_0^2+\mu r_1^2}{2\eta(t)}+O(\eta(t)^{-2})\right)
\end{align*}
for large enough $t$. 
This implies that 
\begin{align*}
-\frac{|x_\mu(t)|^2}{4t}
 & =-\eta(t)\left(1+\frac{(1-\mu)r_0^2+\mu r_1^2}{\eta(t)}+O(\eta(t)^{-2})\right)\\
 & =-\eta(t)-(1-\mu)r_0^2-\mu r_1^2+O(\eta(t)^{-1})
\end{align*}
for large enough $t$. Then 
\begin{equation}
\label{eq:3.18}
\begin{split}
U_d(x_\mu(t),t) & =(4\pi t)^{d-\frac{n}{2}}e^{-\frac{|x_\mu(t)|^2}{4t}}
\int_{{\mathbb R}^n}e^{\frac{2\langle x_\mu(t),y\rangle-|y|^2}{4t}}\phi(y)\,dy\\
 & \to me^{-(1-\mu)r_0^2-\mu r_1^2}\int_{{\mathbb R}^n}\phi(y)\,dy
 =ke^{-(1-\mu)r_0^2-\mu r_1^2}
\end{split}
\end{equation}
for large enough $t$. 
Therefore, by \eqref{eq:3.16}, \eqref{eq:3.17}, and \eqref{eq:3.18} 
we obtain 
$$
F\left(ke^{-(1-\mu)r_0^2-\mu r_1^2}\right)
\ge(1-\mu)F\left(ke^{-r_0^2}\right)+\mu F\left(ke^{-r_1^2}\right).
$$
Since $r_0$, $r_1\in[0,\infty)$ and $\mu\in(0,1)$ are arbitray, 
we deduce that the function
$$
[0,\infty)\ni r\mapsto F\left(ke^{-r^2}\right)
$$
is concave for any $k>0$. 
Therefore we deduce from Lemma~\ref{Lemma:2.5} that $F$-concavity is weaker than log-concavity in ${\mathcal A}([0,\infty))$. 

On the other hand, since $U_d[\phi]$ possesses the eventual log-concavity property for all $\phi\in {\mathcal L}$ 
(see e.g., Theorem~\ref{Theorem:1.1} with $\beta=1$ and \cite{LV}*{Corollary~3.4}), 
$U_d[\phi]$ possesses the eventual $F$-concavity property for all $\phi\in {\mathcal L}$ 
if $F$-concavity is weaker than log-concavity in ${\mathcal A}([0,\infty))$. 
The proof of Theorem~\ref{Theorem:1.5} is complete.
$\Box$
\vspace{5pt}

Similarly, we prove Theorem~\ref{Theorem:1.4}.
\vspace{5pt}
\newline
{\bf Proof of Theorem~\ref{Theorem:1.4}.}
We apply the same argument as in the proof of Theorem~\ref{Theorem:1.5} with $d=n/2$ and $k=M_\phi$. 
Then 
\begin{equation}
\label{eq:3.19}
m=1,\quad \eta(t)=0,\quad x(t:r)=2\sqrt{t}re_1,\quad
\lim_{t\to\infty}U_d(x(t:r),t)=M_\phi e^{-r^2},
\end{equation}
for $t>0$ and $r\in[0,\infty)$. 

Assume that $U_{n/2}[\phi]$ possesses the eventual $F$-concavity property for all $\phi\in{\mathcal L}$. 
Let $r_0$, $r_1\in[0,\infty)$. Set
$$
x_\mu(t):=(1-\mu)x_0(t:r_0)+\mu x_1(t:r_1)
=x(t:(1-\mu)r_0+\mu r_1).
$$
Then, for any $\phi\in{\mathcal L}$, we see that 
$$
F(U_d(x_\mu(t),t))\ge(1-\mu)F(U_d(x_0(t),t))+\mu F(U_d(x_1(t),t))
$$
for large enough $t$. This together with \eqref{eq:3.19} implies that 
$$
F\left(M_\phi e^{-(1-\mu)r_0^2-\mu r_1^2}\right)\ge (1-\mu)F\left(M_\phi e^{-r_0^2}\right)+\mu F\left(M_\phi e^{-r_1^2}\right),
\quad r_0, r_1\in[0,\infty).
$$
Similarly to the proof of Theorem~\ref{Theorem:1.5}, 
applying Lemma~\ref{Lemma:2.5}, 
we see that $F$-concavity is weaker than log-concavity in ${\mathcal A}([0,\infty))$. 
Furthermore, 
we observe from the eventual log-concavity property of $U_d[\phi]$ 
that $U_d[\phi]$ possesses the eventual $F$-concavity property for all $\phi\in{\mathcal L}$ 
if $F$-concavity is weaker than log-concavity in ${\mathcal A}([0,\infty))$. Thus Theorem~\ref{Theorem:1.4} follows.
$\Box$
\section{Eventual $F$-concavity with $d\le n/2$ in ${\mathcal L}_A$}
In this section, under condition~($\mbox{F}_a$) with $a\in(0,\infty)$, 
we study the eventual $F$-concavity property of $U_d[\phi]$ with $\phi\in{\mathcal L}_A$ 
to prove Theorem~\ref{Theorem:1.6}. 
Similarly to Section~3, we can assume, without loss of generality, that a = 1.
We use the same notation as in Section~3. 
\subsection{Proof of Theorem~\ref{Theorem:1.6}~(1)}
We first prove Theorem~\ref{Theorem:1.6}~(1) using Lemma~\ref{Lemma:2.4}. 
\vspace{5pt}
\newline
{\bf Proof Theorem~\ref{Theorem:1.6}~(1).} 
Assume condition~($\mbox{F}_a$) with $a=1$. Furthermore, assume that the $F$-concavity property is preserved by the heat flow in ${\mathbb R}^n$. 
Then, 
by Lemma~\ref{Lemma:2.4} we see that, 
for any $k\in(0,1)$, the function $[0,\infty)\ni r\mapsto F(k e^{-r^2})$ is concave, that is, 
$$
\frac{d^2}{dr^2}F(k e^{-r^2})\le 0,\quad r\in(0,\infty). 
$$
Since $F(k e^{-r^2})={\mathcal F}(-\log k+r^2)$ for $r\in(0,\infty)$, 
we have
$$
\frac{d^2}{dr^2}F(k e^{-r^2})=4{\mathcal F}''(-\log k+r^2)r^2+2{\mathcal F}'(-\log k+r^2)\le 0
$$
for $r\in(0,\infty)$. Since ${\mathcal F}'(r)<0$ for $r\in(0,\infty)$, we see that 
$$
-\frac{1}{2}\le \frac{r^2{\mathcal F}''(-\log k+r^2)}{{\mathcal F}'(-\log k+r^2)}=\frac{(s+\log k){\mathcal F}''(s)}{{\mathcal F}'(s)}
=-\left(1+\frac{\log k}{s}\right)\kappa_F(s)
$$
for $r\in(0,\infty)$, where $s:=-\log k+r^2$. 
This implies that
$$
\kappa_F(s)\le\frac{1}{2}\left(1+\frac{\log k}{s}\right)^{-1}=\frac{1}{2}-\frac{1}{2}\frac{\log k}{s}+O\left(s^{-2}\right)
$$
for large enough $s$. Then 
$$
\sigma_F(s)=s\left(\kappa_F(s)-\frac{1}{2}\right)\le -\frac{\log k}{2}+O(s^{-1})
$$
for large enough $s$. 
This implies the desired conclusion. The proof is complete. 
$\Box$
%
\subsection{Eventual $F$-concavity for ${\mathcal L}_A$}
We assume condition~($\mbox{F}_a$) with $a=1$ and $\kappa_F^*=1/2$. 
Let $\phi\in{\mathcal L}_A$ and $d\le n/2$.
Then, for any $\epsilon>0$, there exists $R_\epsilon>0$ such that 
$$
\kappa_F(r)-\frac{1}{2}<\epsilon
$$
for $r\in(R_\epsilon,\infty)$. 
Let $\kappa_*\in(0, 1/2)$. Taking small enough $\epsilon$ if necessary, we have
$$
\nu_F(r)=\frac{r}{\kappa_F(r)}\left(\kappa_F(r)-\frac{1}{2}\right)\le\frac{\epsilon}{\kappa_*}r\le\frac{1}{2}r
$$
for $r\in(R_\epsilon,\infty)$ with $\kappa_F(r)\ge\kappa_*$. 
Then we find $T>0$ and $m>0$ such that 
$$
v_d(x,t)>\nu_F(v_d(x,t))=v_d(x,t)\left(1-\frac{1}{2\kappa_F(v_d(x,t))}\right)
$$
for $(x,t)\in {\mathbb R}^n\times[T,\infty)$ with $\kappa_F(v_d(x,t))\ge\kappa_*$ and $v_d(x,t)\ge m$. 
Furthermore, 
$$
\frac{{\mathcal F}''(v_d(x,t))}{{\mathcal F}'(v_d(x,t))}(v_d(x,t)-\nu_F(v_d(x,t)))
=-\kappa_F(v_d(x,t))\left(1-\left(1-\frac{1}{2\kappa_F(v_d(x,t))}\right)\right)=-\frac{1}{2}
$$
for $(x,t)\in {\mathbb R}^n\times[T,\infty)$ with $\kappa_F(v_d(x,t))\ge\kappa_*$ and $v_d(x,t)\ge m$. 
These imply that 
\begin{equation*}
\begin{split}
\frac{\partial^2}{\partial\xi^2}F\left(U_d[\phi](x,t)\right)
 & ={\mathcal F}''(v_d(x,t))\left(\frac{\partial}{\partial\xi} v_d(x,t)\right)^2+{\mathcal F}'(v_d(x,t))\frac{\partial^2}{\partial\xi^2}v_d(x,t)\\
 & =\frac{{\mathcal F}'(v_d(x,t))}{v_d(x,t)-\nu_F(v_d(x,t))}J_F[v_d](x,t;\xi)
\end{split}
\end{equation*}
for $(x,t)\in {\mathbb R}^n\times[T,\infty)$ with $\kappa_F(v_d(x,t))\ge\kappa_*$ and $v_d(x,t)\ge m$,  
where 
\begin{equation}
\label{eq:4.1}
\begin{split}
J_F[v_d](x,t;\xi)
 & :=-\frac{1}{2}\left(\frac{\partial}{\partial\xi} v_d(x,t)\right)^2+(v_d(x,t)-\nu_F(v_d(x,t)))\frac{\partial^2}{\partial\xi^2}v_d(x,t)\\
 & =-\frac{1}{2}\left(\frac{\langle x,\xi\rangle}{2t}-\left(\frac{\partial}{\partial\xi}\log w\right)(x,t)\right)^2\\
 & \qquad\quad
 +\left(\left(\frac{n}{2}-d\right)\log(4\pi t)+\frac{|x|^2}{4t}-\log w(x,t)-\nu_F(v_d(x,t))\right)\\
 & \qquad\qquad\qquad
 \times\left(\frac{1}{2t}-\left(\frac{\partial^2}{\partial\xi^2}\log w\right)(x,t)\right).
\end{split}
\end{equation}
Therefore, in the case of $\kappa_F^*=1/2$, we observe that 
\begin{equation}
\label{eq:4.2}
\frac{\partial^2}{\partial\xi^2}F\left(U_d[\phi](x,t)\right)
\le 0\quad\mbox{if and only if}\quad J_F[v_d](x,t;\xi)\ge 0
\end{equation}
for $(x,t)\in {\mathbb R}^n\times[T,\infty)$ with $\kappa_F(v_d(x,t))\ge\kappa_*$ and $v_d(x,t)\ge m$. 
The main purpose of this subsection is to prove the following proposition. 
\begin{proposition}
\label{Proposition:4.1}
Assume condition~{\rm ($\mbox{F}_a$)} with $a=1$ and $\kappa_F^*=1/2$. 
Let $\phi\in{\mathcal L}_A$ and $\kappa_*\in(0,1/2)$. 
Consider the following three cases: 
\begin{itemize}
  \item[{\rm (C1)}] 
  $\displaystyle{d<\frac{n-1}{2}}$ and 
  \begin{equation}
  \label{eq:4.3}
  \limsup_{r\to\infty}\,\left(\nu_F(r)-\frac{1}{2}\log r\right)<\infty;
  \end{equation}
  \item[{\rm (C2)}] 
  $\displaystyle{\frac{n-1}{2}\le d<\frac{n}{2}}$ and 
  \begin{equation}
  \label{eq:4.4}
  \lim_{r\to\infty}\,\left(\nu_F(r)-\left(\frac{n}{2}-d\right)\log r\right)=-\infty;
  \end{equation}
  \item[{\rm (C3)}] 
  $\displaystyle{d=\frac{n}{2}}$, $M_\phi<1$, 
  \begin{equation}
  \label{eq:4.5}
  \kappa_F(r)\le\frac{1}{2}\quad\mbox{for $r>0$},
  \quad\mbox{and}\quad
  \lim_{r\to\infty}\,\nu_F(r)=-\infty.
  \end{equation}
\end{itemize} 
Then there exist $T>0$ and $m>0$ such that 
\begin{equation}
\label{eq:4.6}
\frac{\partial^2}{\partial \xi^2}F\left(U_d[\phi](x,t)\right)\le 0,\quad \xi\in{\mathbb S}^{n-1},
\end{equation}
for $(x,t)\in {\mathbb R}^n\times[T,\infty)$ with $\kappa_F(v_d(x,t))\ge\kappa_*$ and $v_d(x,t)\ge m$.
\end{proposition}
We remark that 
conditions~\eqref{eq:4.3}, \eqref{eq:4.4}, and \eqref{eq:4.5} are equivalent to \eqref{eq:1.8}, \eqref{eq:1.9}, and \eqref{eq:1.10}, respectively 
(see Lemma~\ref{Lemma:2.7}). 
\begin{lemma}
\label{Lemma:4.1}
Let $\phi\in{\mathcal L}$. 
Assume condition~{\rm ($\mbox{F}_a$)} with $a=1$. 
Consider either {\rm (C1)}, {\rm (C2)}, or {\rm (C3)}. 
Then, for any $R>0$, 
there exist $T>0$ and $m>0$ such that 
\begin{equation}
\label{eq:4.7}
J_F[v_d](x,t;\xi)\ge 0,\quad \xi\in{\mathbb S}^{n-1},
\end{equation}
for $x\in B(0,Rt)$ and $t\in[T,\infty)$ with $v_d(x,t)\ge m$. 
\end{lemma}
{\bf Proof.}
Let $\phi\in{\mathcal L}$, $\xi\in{\mathbb S}^{n-1}$, and $R>0$. 
Consider either (C1), (C2), or (C3). 
Let $M$ be large enough. 
Then we find $T>0$ and $m>0$ such that 
\begin{equation}
\label{eq:4.8}
\psi_d(x,t)\ge M
\end{equation}
for $(x,t)\in{\mathbb R}^n\times[T,\infty)$ with $v_d(x,t)\ge m$, where
\begin{equation}
\label{eq:4.9}
\psi_d(x,t):=
\left\{
\begin{array}{ll}
\displaystyle{\left(\frac{n-1}{2}-d\right)\log t} & \quad\mbox{if}\quad\displaystyle{d<\frac{n-1}{2}},\vspace{7pt}\\
\displaystyle{\left(\frac{n}{2}-d\right)\log v_d(x,t)-\nu_F(v_d(x,t))} & \quad\mbox{if}\quad\displaystyle{\frac{n-1}{2}\le d\le\frac{n}{2}},\vspace{3pt}
\end{array}
\right.
\end{equation}
for $(x,t)\in{\mathbb R}^n\times(0,\infty)$. 
\vspace{3pt}
\newline
\underline{Step 1.} 
Let $T>0$ and $m>0$ be large enough. 
We show that 
\begin{equation}
\label{eq:4.10}
\nu_F(v_d(x,t))
\le\left(\frac{n}{2}-d\right)\log t-\psi_d(x,t)+C
\end{equation}
for $x\in B(0,Rt)$ and $t\in[T,\infty)$ with $v_d(x,t)\ge m$. 
In the case of (C1), 
by Lemma~\ref{Lemma:2.8}, \eqref{eq:1.11}, and \eqref{eq:4.9} we obtain 
\begin{equation*}
\begin{split}
\nu_F(v_d(x,t)) & \le\frac{1}{2}\log v_d(x,t)+C\\
 & =\frac{1}{2}\log\left(\left(\frac{n}{2}-d\right)\log(4\pi t)+\frac{|x|^2}{4t}-\log w(x,t)\right)+C\\
 & \le \frac{1}{2}\log(Ct)+C\\
 & \le \left(\frac{n}{2}-d\right)\log t-\psi_d(x,t)+C
\end{split}
\end{equation*}
for $x\in B(0,Rt)$ and $t\in[T,\infty)$ with $v_d(x,t)\ge m$. 
Thus \eqref{eq:4.10} holds in the case of (C1). 
In the case of (C2), 
by Lemma~\ref{Lemma:2.8}, \eqref{eq:1.11}, and \eqref{eq:4.9} we have
\begin{equation*}
\begin{split}
\nu_F(v_d(x,t)) & =\left(\frac{n}{2}-d\right)\log\left(\left(\frac{n}{2}-d\right)\log(4\pi t)+\frac{|x|^2}{4t}-\log w(x,t)\right)-\psi_d(x,t)\\
 & \le \left(\frac{n}{2}-d\right)\log(Ct)-\psi_d(x,t)
 =\left(\frac{n}{2}-d\right)\log t-\psi_d(x,t)+C
 \end{split}
\end{equation*}
for $x\in B(0,Rt)$ and $t\in[T,\infty)$.  
Thus \eqref{eq:4.10} holds in the case of (C2). 
In the case of (C3), by \eqref{eq:4.9} we easily obtain \eqref{eq:4.10}. 
Thus \eqref{eq:4.10} holds. 
\vspace{3pt}
\newline
\underline{Step 2.} 
We complete the proof of Lemma~\ref{Lemma:4.1}. 
Let $T>0$, $M>0$, and $m>0$ be large enough such that \eqref{eq:4.8} and \eqref{eq:4.10} hold. 
Let $\xi\in{\mathbb S}^{n-1}$. 
Taking large enough $M$ if necessary, 
by Lemma~\ref{Lemma:2.8} and \eqref{eq:4.1} we have
\begin{equation*}
\begin{split}
 & J_F[v_d](x,t;\xi)  
=-\frac{1}{2}\left(\frac{\langle x,\xi\rangle}{2t}+O(t^{-1})\right)^2\\
 & \qquad
 +\left(\left(\frac{n}{2}-d\right)\log(4\pi t)+\frac{|x|^2}{4t}+O(t^{-1}|x|)+O(1)-\nu_F(v_d(x,t))\right)\left(\frac{1}{2t}+O(t^{-2})\right)\\
 & \ge-\frac{1}{2}\left(\frac{\langle x,\xi\rangle}{2t}+O(t^{-1})\right)^2\\
 & \qquad
 +\left(\psi(v_d(x,t))+O(1)+\frac{|x|^2}{4t}+O(t^{-1}|x|))\right)\left(\frac{1}{2t}+O(t^{-2})\right)\\
 & \ge\frac{1}{2t}\psi(v_d(x,t))(1+O(t^{-1}))-\frac{\langle x,\xi\rangle^2}{8t^2}+\frac{|x|^2}{8t^2}
 +O(t^{-1})+O(t^{-2}|x|)+O(t^{-3}|x|^2)\\
 & \ge\frac{M}{4t}-Ct^{-1}-C(R+R^2) t^{-1}>0
 \end{split}
\end{equation*}
for $x\in B(0,Rt)$ and $t\in[T,\infty)$ with $v_d(x,t)\ge m$. 
This implies \eqref{eq:4.7}. The proof is complete. 
$\Box$\vspace{5pt}
\begin{lemma}
\label{Lemma:4.2}
Let $\phi\in{\mathcal L}_A$.
Assume condition~{\rm ($\mbox{F}_a$)} with $a=1$. 
Consider either {\rm (C1)}, {\rm (C2)}, or {\rm (C3)}. 
Then, for any large enough $R>0$, 
there exist $T>0$ and $m>0$ such that 
\begin{equation}
\label{eq:4.11}
J_F[v_d](x,t;\xi)\ge 0,\quad \xi\in{\mathbb S}^{n-1},
\end{equation}
for $x\in {\mathbb R}^n\setminus B(0,Rt)$ and $t\in[T,\infty)$ with $v_d(x,t)\ge m$.
\end{lemma}
{\bf Proof.}
Let $\phi\in{\mathcal L}_A$ and $\xi\in{\mathbb S}^{n-1}$. 
Let $R>0$ be large enough. 
We can assume, without loss of generality, that
\begin{equation}
\label{eq:4.12}
\begin{split}
 & P_\phi\subset(-\infty,0)\times{\mathbb R}^{n-1},\quad
 \partial P_\phi\cap (\{0\}\times{\mathbb R}^{n-1})\not=\emptyset,\quad \mbox{$x=(X,0)\in{\mathbb R}^n$ with $X\ge 0$},\\
 & \xi=e_1,\quad \frac{\partial}{\partial\xi}=\frac{\partial}{\partial x_1},\quad\mbox{if $n=1$},\\
 & 
\xi=(\cos\theta) e_1+(\sin\theta) e_2,
\quad
\frac{\partial}{\partial\xi}=\cos\theta\frac{\partial}{\partial x_1}+\sin\theta\frac{\partial}{\partial x_2},\quad\mbox{if $n\ge 2$}, 
\end{split}
\end{equation}
where $\theta\in(-\pi/2,\pi/2]$. We set $\theta=0$ if $n=1$, for the convenience. 

Let $T$ be large enough. 
By Lemmas~\ref{Lemma:2.8} and \ref{Lemma:2.9} we have
\begin{equation}
\label{eq:4.13}
\log w(x,t)\le -\log\frac{X}{t}+C
\end{equation}
for $x=(X,0)$ with $X>0$ and $t\in[T,\infty)$ and 
\begin{equation}
\label{eq:4.14}
\begin{split}
 & \left|\left(\frac{\partial}{\partial\xi}\log w\right)(x,t)\right|\le CX^{-1}\cos\theta+Ct^{-1}|\sin\theta|,\\
 & \left|\left(\frac{\partial^2}{\partial\xi^2}\log w\right)(x,t)\right|\le CX^{-2}\cos^2\theta+Ct^{-2}\sin^2\theta,\\
\end{split}
\end{equation}
for $x=(X,0)$ with $X\ge Rt$ and $t\in[T,\infty)$. 
\vspace{3pt}
\newline
\underline{Step 1.} 
We find $m>0$ such that 
\begin{equation}
\label{eq:4.15}
\left(\frac{n}{2}-d\right)\log(4\pi t)-\log w(x,t)-\nu_F(v_d(x,t))
\ge\psi(v_d(x,t))-C
\end{equation}
for $x=(X,0)$ and $t\in[T,\infty)$ with $X\ge Rt$ and $v_d(x,t)\ge m$.
In the case of (C1),
by \eqref{eq:2.6}, \eqref{eq:4.9}, and \eqref{eq:4.13} we find $m>0$ such that 
\begin{equation*}
\begin{split}
 & \left(\frac{n}{2}-d\right)\log(4\pi t)-\log w(x,t)-\nu_F(v_d(x,t))\\
 & \ge \left(\frac{n}{2}-d\right)\log(4\pi t)+\log\frac{X}{t}-C\\
 & \qquad\quad
 -\frac{1}{2}\log\left(\left(\frac{n}{2}-d\right)\log(4\pi t)+\frac{X^2}{4t}-\log w(x,t)\right)-C\\
 & \ge \left(\frac{n}{2}-d\right)\log(4\pi t)+\log\frac{X}{t}-\frac{1}{2}\log\frac{X^2}{t}-C\\
 & =\left(\frac{n-1}{2}-d\right)\log t-C
 =\psi(v_d(x,t))-C
\end{split}
\end{equation*}
for $x=(X,0)$ and $t\in[T,\infty)$ with $X\ge Rt$ and $v_d(x,t)\ge m$.
This implies \eqref{eq:4.15} in the case of (C1). 
In the case of (C2), 
by \eqref{eq:2.6}, \eqref{eq:4.9}, and \eqref{eq:4.13} we have 
\begin{align*}
 & \left(\frac{n}{2}-d\right)\log(4\pi t)-\log w(x,t)-\nu_F(v_d(x,t))\\
 & \ge \left(\frac{n}{2}-d\right)\log(4\pi t)+\log\frac{X}{t}-C\\
 & \qquad\quad
 -\left(\frac{n}{2}-d\right)\log\left(\left(\frac{n}{2}-d\right)\log(4\pi t)+\frac{X^2}{4t}-\log w(x,t)\right)+\psi_d(x,t)\\
 & \ge \left(\frac{n}{2}-d\right)\log(4\pi t)+\log\frac{X}{t}-\left(\frac{n}{2}-d\right)\log\frac{X^2}{t}-C+\psi_d(x,t)\\
 & =\left(1-2\left(\frac{n}{2}-d\right)\right)\log\frac{X}{t}+\psi_d(x,t)-C
 \ge \psi_d(x,t)-C
\end{align*}
for $x=(X,0)$ with $X\ge Rt$ and $t\in[T,\infty)$. 
This implies \eqref{eq:4.15} in the case of (C2). 
In the case of (C3), 
by \eqref{eq:4.13} we see that
\begin{align*}
 & \left(\frac{n}{2}-d\right)\log(4\pi t)-\log w(x,t)-\nu_F(v_d(x,t))\\
 & \ge \log\frac{X}{t}-C+\psi_d(x,t)\ge \psi_d(x,t)-C
\end{align*}
for $x=(X,0)$ with $X\ge Rt$ and $t\in[T,\infty)$.
This implies \eqref{eq:4.15} in the case of (C3). Thus  \eqref{eq:4.15} holds. 
\vspace{3pt}
\newline
\underline{Step 2.} 
We complete the proof of Lemma~\ref{Lemma:4.2}. 
Let $M$ be large enough. 
Let $T>0$ and $m>0$ be large enough such that \eqref{eq:4.8} and \eqref{eq:4.15} hold. 
By \eqref{eq:1.11}, \eqref{eq:4.14}, and \eqref{eq:4.15} 
we obtain 
\begin{equation}
\label{eq:4.16}
\begin{split}
 & (v_d(x,t)-\nu_F(v_d(x,t)))\frac{\partial^2}{\partial\xi^2}v_d(x,t)\\
 & \ge\left(\frac{X^2}{4t}+\psi(v_d(x,t))-C\right)\left(\frac{1}{2t}+O(X^{-2}\cos^2\theta)+O(t^{-2}\sin^2\theta)\right)\\
 & \ge \frac{X^2}{8t^2}+O(t^{-1})+O(t^{-3}X^2\sin^2\theta)\\
 & \qquad\quad
 +\frac{1}{2t}\psi_d(x,t)\left(1+O(tX^{-2}\cos^2\theta)+O(t^{-1}\sin^2\theta)\right)\\
 & \ge \frac{X^2}{8t^2}+\frac{1}{4t}\psi_d(x,t)+O(t^{-1})+O(t^{-3}X^2\sin^2\theta)
\end{split}
\end{equation}
for $x=(X,0)\in{\mathbb R}^n$ and $t\in[T,\infty)$ with $X\ge Rt$ and $v_d(x,t)\ge m$. 
Similarly, by \eqref{eq:4.14} we have 
\begin{equation}
\label{eq:4.17}
\begin{split}
 & \left(\frac{\partial}{\partial\xi} v_d(x,t)\right)^2\\
 & =\left(\frac{X\cos\theta}{2t}+O(X^{-1}\cos\theta)+O(t^{-1}|\sin\theta|)\right)^2\\
 & =\frac{X^2\cos^2\theta}{4t^2}+O(X^{-2})+O(t^{-2}\sin^2\theta)+O(t^{-1})+O(t^{-1}X^{-1})+O(t^{-2}X|\sin\theta|)\\
 & =\frac{X^2\cos^2\theta}{4t^2}+O(t^{-1})+O(t^{-3}X^2\sin^2\theta)
\end{split}
\end{equation}
for $x=(X,0)\in{\mathbb R}^n$ and large enough $t$ with $X\ge Rt$. 
Therefore, by \eqref{eq:4.1}, \eqref{eq:4.8}, \eqref{eq:4.16}, and \eqref{eq:4.17}, 
taking a large enough $T$ if necessary, we obtain 
\begin{equation*}
\begin{split}
J_F[v_d](Xe_1,t;\xi) & \ge \frac{X^2}{8t^2}\sin^2\theta+\frac{1}{4t}\psi_d(Xe_1,t)+O(t^{-1})+O(t^{-3}X^2\sin^2\theta)\\
 & \ge\frac{X^2}{16t^2}\sin^2\theta+\frac{1}{4t}M-Ct^{-1}
\end{split}
\end{equation*}
for $X\ge Rt$ and $t\in[T,\infty)$ with $v_d(x,t)\ge m$. 
Taking large enough $M$ if necessary, we obtain \eqref{eq:4.11}. 
Thus Lemma~\ref{Lemma:4.2} follows.
$\Box$
\vspace{3pt}
\newline
{\bf Proof of Proposition~\ref{Proposition:4.1}.}
Consider either (C1), (C2), or (C3). 
By Lemmas~\ref{Lemma:4.1} and \ref{Lemma:4.2} we find $T>0$ and $m>0$ such that 
$$
J_F[v_d](x,t;\xi)\ge 0,\quad \xi\in{\mathbb S}^{n-1},
$$
for $(x,t)\in {\mathbb R}^n\times[T,\infty)$ with $v_d(x,t)\ge m$. 
This together with \eqref{eq:4.2} implies \eqref{eq:4.6}. 
Thus Proposition~\ref{Proposition:4.1} follows.
$\Box$\vspace{5pt}
\newline
Then we obtain 
\begin{proposition}
\label{Proposition:4.2}
Assume condition~{\rm ($\mbox{F}_a$)} with $a=1$. 
Let $\phi\in{\mathcal L}_A$. 
Consider either {\rm (C1)}, {\rm (C2)}, or {\rm (C3)}. 
Then $U_d[\phi]$ possesses the eventual $F$-concavity property, that is, 
there exists $T>0$ such that 
\begin{equation}
\label{eq:4.18}
\frac{\partial^2}{\partial \xi^2}F\left(U_d[\phi](x,t)\right)\le 0,
\quad \xi\in{\mathbb S}^{n-1},
\end{equation}
for $(x,t)\in {\mathbb R}^n\times[T,\infty)$. 
\end{proposition}
{\bf Proof.}
Consider either (C1), (C2), or (C3). Then $\kappa_F^*\le 1/2$.  
It suffices to consider the case of $\kappa_F^*=1/2$. 
Indeed, 
if $\kappa_F^*<1/2$, 
then Proposition~\ref{Proposition:4.2} follows from Theorem~\ref{Theorem:1.2}~(1) and Theorem~\ref{Theorem:1.3}~(1). 

Let $\kappa_F^*=1/2$ and $\kappa_*\in(0,1/2)$. 
In cases (C1) and (C2), 
by Lemma~\ref{Lemma:3.1}~(1) and Proposition~\ref{Proposition:4.1}
we find $T>0$ such that 
\eqref{eq:4.18} holds for $(x,t)\in {\mathbb R}^n\times[T,\infty)$ with $\kappa_F(v_d(x,t))\ge\kappa_*$.
This together with Proposition~\ref{Proposition:3.1}~(2) implies that \eqref{eq:4.18} holds 
for $(x,t)\in {\mathbb R}^n\times[T,\infty)$. 
Thus Proposition~\ref{Proposition:4.2} holds in cases (C1) and (C2). 

In case (C3), by Lemma~\ref{Lemma:3.1}~(2) and Proposition~\ref{Proposition:4.1}, 
we find $T'>0$ and $L>0$ such that 
\eqref{eq:4.18} holds for $x\in{\mathbb R}^n\setminus B(0,L\sqrt{t})$ and $t\in[T,\infty)$ with $\kappa_F(v_d(x,t))\ge\kappa_*$.
This together with Proposition~\ref{Proposition:3.1}~(1) and (2) imply that \eqref{eq:4.18} holds 
for $(x,t)\in {\mathbb R}^n\times[T,\infty)$. 
Thus Proposition~\ref{Proposition:4.2} holds in case (C3).
The proof is complete.
$\Box$
\subsection{Proofs of Theorem~\ref{Theorem:1.6}~(2), (3), and (4)}
{\bf Proof of Theorem~\ref{Theorem:1.6}~(2).}
Let $d<(n-1)/2$ and $\phi\in{\mathcal L}_A$. 
By Remark~\ref{Remark:1.1}~(1) it suffices to consider the case of $a=1$. 
\vspace{3pt}
\newline
\underline{Step 1}.
Assume that \eqref{eq:1.8} holds. 
It follows from Lemma~\ref{Lemma:2.7} that \eqref{eq:4.3} holds. 
Proposition~\ref{Proposition:4.2} implies that $U_d[\phi]$ possesses the eventual $F$-concavity property. 
\vspace{3pt}
\newline
\underline{Step 2}.
Assume that \eqref{eq:1.8} does not hold. 
Then $\kappa_F^*\ge 1/2$. 
Let $\ell$ be large enough, and set 
$$
\phi_\ell:=\ell\chi_\Omega,\quad\mbox{where}\quad \Omega:=\left(-\frac{1}{2\ell},0\right)\times(-1,0)^{n-1}.
$$ 
Then Lemma~\ref{Lemma:2.10} implies that $\phi_\ell\in{\mathcal L}_A$. 
We prove that $U_d[\phi_\ell](t)$ does not possess the eventual $F$-concavity property. 
By Theorem~\ref{Theorem:1.2}~(2) it suffices to consider the case of $\kappa_F^*=1/2$.

Let $X>0$, and set $X_t:=X/2t$. 
It follows that
\begin{equation}
\label{eq:4.19}
\begin{split}
w(Xe_1,t) & =\ell\int_\Omega e^{\frac{2Xy_1-|y|^2}{4t}}\,dy
=\ell(1+O(t^{-1}))\int_{-1/2\ell}^0 e^{X_ty_1}\,dy_1\\
 & =\frac{\ell}{X_t}(1+O(t^{-1}))\left(1-e^{-\frac{X_t}{2\ell}}\right)
\end{split}
\end{equation}
for large enough $t$. 
Taking large enough $L>0$, we have
\begin{equation}
\label{eq:4.20}
-\log w(Xe_1,t)=\log X_t-\log\ell(1+o(1))\,\,
\left\{
\begin{array}{l}
\le \displaystyle{\log X_t-\frac{1}{2}\log\ell},\vspace{5pt}\\
\ge \displaystyle{\log X_t-2\log\ell},
\end{array}
\right. 
\end{equation}
for large enough $t$ and $X_t\ge L$.
Since
\begin{align*}
\left(\frac{\partial}{\partial x_1}w\right)(Xe_1,t) & =\frac{\ell}{2t}(1+O(t^{-1}))\int_{-1/2\ell}^0 y_1e^{X_ty_1}\,dy_1\\
 & =\frac{\ell}{2t}(1+O(t^{-1}))\left[\frac{1}{2\ell X_t}e^{-\frac{X_t}{2\ell}}-\frac{1}{X_t}\int_{-1/2\ell}^0e^{X_ty_1}\,dy_1\right]\\
 & =\frac{\ell}{2t}(1+O(t^{-1}))\left[\frac{1}{2\ell X_t}e^{-\frac{X_t}{2\ell}}-\frac{1}{X_t^2}\left(1-e^{-\frac{X_t}{2\ell}}\right)\right]\\\
 & =-\frac{\ell}{2tX_t^2}(1+O(t^{-1}))(1+O(X_te^{-\frac{X_t}{2\ell}}))
\end{align*}
for $X_t\ge L$ and large enough $t$, by \eqref{eq:4.19} we have
$$
\left(\frac{\partial}{\partial x_1}\log w\right)(Xe_1,t)=-\frac{1}{2 tX_t}(1+O(X_te^{-\frac{X_t}{2\ell}})(1+O(t^{-1}))
$$
for $X_t\ge L$ and large enough $t$.
This together with \eqref{eq:1.11} implies that
\begin{align*}
\left(\frac{\partial}{\partial x_1} v_d(x,t)\right)^2
 & =\left(X_t+\frac{1}{2 tX_t}(1+O(X_te^{-\frac{X_t}{2\ell}}))(1+O(t^{-1}))\right)^2\\
 & =X_t^2+\frac{1}{4 t^2X_t^2}(1+o(1))
+\frac{1}{t}(1+o(1))
=\frac{X^2}{4t^2}+\frac{1+o(1)}{X^2}+\frac{1+o(1)}{t}
\end{align*}
for $X_t\ge L$ and large enough $t$.
Then, by Lemma~\ref{Lemma:2.9}, \eqref{eq:4.1}, and \eqref{eq:4.20} we have
\begin{equation}
\label{eq:4.21}
\begin{split}
 & J_F[v_d](Xe_1,t;e_1)\\
 & \equiv-\frac{1}{2}\left(\frac{\partial}{\partial x_1} v_d(x,t)\right)^2+(v_d(x,t)-\nu_F(v_d(x,t)))\frac{\partial^2}{\partial x_1^2}v_d(x,t)\\
 & \le -\frac{X^2}{8t^2}-\frac{1+o(1)}{2X^2}-\frac{1+o(1)}{2t}+\frac{X^2}{4t}\left(\frac{1}{2t}-\left(\frac{\partial^2}{\partial x_1^2}\log w\right)(x,t)\right)\\
 & \qquad
 +\left(\left(\frac{n}{2}-d\right)\log(4\pi t)+\log\frac{X}{2t}-\frac{1}{2}\log\ell-\nu_F(v_d(x,t))\right)\left(\frac{1}{2t}+O(X^{-2})\right)\\
  & \le\left(\left(\frac{n}{2}-d\right)\log(4\pi t)+\log\frac{X}{2t}-\frac{1}{2}\log\ell-\nu_F(v_d(x,t))\right)\left(\frac{1}{2t}+O(X^{-2})\right)
\end{split}
\end{equation}
for $X_t\ge L$ and large enough $t$.
Furthermore, 
setting $\tilde{\nu}_F(r)=\nu_F(r)-(\log r)/2$, 
by \eqref{eq:4.20} we have
\begin{align*}
 & \left(\frac{n}{2}-d\right)\log(4\pi t)+\log\frac{X}{2t}-\frac{1}{2}\log\ell-\nu_F(v_d(Xe_1,t))\\
 & \le \left(\frac{n}{2}-d\right)\log(4\pi t)+\log\frac{X}{2t}-\frac{1}{2}\log\ell\\
 & \qquad\quad
 -\frac{1}{2}\log\left[\left(\frac{n}{2}-d\right)\log(4\pi t)+\frac{X^2}{4t}+\log\frac{X}{2t}-2\log\ell\right]-\tilde{\nu}_F(v_d(Xe_1,t))\\
 & \le \left(\frac{n}{2}-d\right)\log(4\pi t)+\log\frac{X}{2t}-\frac{1}{2}\log\frac{X^2}{4t}+C-\tilde{\nu}_F(v_d(Xe_1,t))\\
 & \le 2\left(\frac{n-1}{2}-d\right)\log t-\tilde{\nu}_F(v_d(Xe_1,t))
\end{align*}
for $X_t\ge L$ and large enough $t$. 
On the other hand, since \eqref{eq:1.8} does not hold, 
by Lemma~\ref{Lemma:2.7} we find $\{r_j\}\subset(0,\infty)$ with 
$\lim_{j\to\infty}r_j=\infty$ such that
$$
\lim_{j\to\infty}\tilde{\nu}(r_j)=\lim_{j\to\infty}\left(\nu_F(r_j)-\frac{1}{2}\log r_j\right)=\infty.
$$ 
Then it follows from $\kappa_F^*=1/2$ that 
\begin{equation}
\label{eq:4.22}
\lim_{j\to\infty} \kappa(r_j)=\frac{1}{2}. 
\end{equation}
By \eqref{eq:1.13a}, applying the intermediate value theorem, 
for any large enough $t>0$, 
we find $X(t)\in(2Lt,\infty)$ and $j=1,2,\dots$ such that 
\begin{equation}
\label{eq:4.23}
v_d(X(t)e_1,t)=r_j,\qquad
\tilde{\nu}_F(r_j)\ge \left[2\left(\frac{n-1}{2}-d\right)+1\right]\log t.
\end{equation}
This together with \eqref{eq:4.21} implies that 
$$
J_F[v_d](X(t)e_1,t;e_1)
\le-\log t\cdot \frac{1}{2t}(1+o(1))<0
$$
for large enough $t$. 
This together with \eqref{eq:4.2} and \eqref{eq:4.22} implies that 
$U_d[\phi_\ell](t)$ does not possess the eventual $F$-concavity property. 
Thus Theorem~\ref{Theorem:1.6}~(2) follows. 
$\Box$\vspace{5pt}
\newline
{\bf Proof of Theorem~\ref{Theorem:1.6}~(3).}
Let $(n-1)/2\le d<n/2$ and $\phi\in{\mathcal L}_A$. 
It suffices to consider the case of $a=1$. 
\vspace{3pt}
\newline
\underline{Step 1}.
Assume that \eqref{eq:1.9} holds. 
It follows from Lemma~\ref{Lemma:2.7} that \eqref{eq:4.4} holds. 
Proposition~\ref{Proposition:4.2} implies that $U_d[\phi]$ possesses the eventual $F$-concavity property. 
\vspace{3pt}
\newline
\underline{Step 2}.
Assume that \eqref{eq:1.9} does not hold. Then $\kappa_F^*\ge 1/2$. 
For any large enough $\ell>0$, 
let $\phi_\ell$ be as in Step 2 in the proof of Theorem~\ref{Theorem:1.6}~(2). 
Then $\phi_\ell\in{\mathcal L}_A$. 
We prove that $U_d[\phi_\ell]$ does not possess the eventual $F$-concavity property. 
By Theorem~\ref{Theorem:1.2}~(2) it suffices to consider the case of $\kappa_F^*=1/2$. 
 
Since \eqref{eq:1.9} does not hold, 
by Lemma~\ref{Lemma:2.7} 
we find $\{r_j\}\subset(0,\infty)$ with 
$\lim_{j\to\infty}r_j=\infty$ such that 
\begin{equation}
\label{eq:4.24}
\lim_{j\to\infty}\left[\left(\frac{n}{2}-d\right)\log r_j-\nu_F(r_j)\right]<\infty.
\end{equation}
Then it follows from $\kappa_F^*=1/2$ that 
\begin{equation}
\label{eq:4.25}
\lim_{j\to\infty}\kappa_F(r_j)=\frac{1}{2}.
\end{equation}
Let $L$ be large enough. 
Taking large enough $\ell$ if necessary, by \eqref{eq:4.9} and \eqref{eq:4.20} we have
\begin{equation}
\label{eq:4.26}
\begin{split}
 & \left(\frac{n}{2}-d\right)\log(4\pi t)-\log w(Xe_1,t)-\nu_F(v_d(Xe_1,t))\\
 & \le \left(\frac{n}{2}-d\right)\log(4\pi t)+\log\frac{X}{2t}-\frac{1}{2}\log\ell\\
 & \qquad\quad
 -\left(\frac{n}{2}-d\right)\log\left[\left(\frac{n}{2}-d\right)\log(4\pi t)+\frac{X^2}{4t}+\log\frac{X}{2t}+O(\ell)\right]+\psi_d(x,t)\\
 & \le \left(\frac{n}{2}-d\right)\log(4\pi t)+\log\frac{X}{2t}-\frac{1}{2}\log\ell-\left(\frac{n}{2}-d\right)\log\frac{X^2}{4t}+C+\psi_d(x,t)\\
 & =\left(1-2\left(\frac{n}{2}-d\right)\right)\log\frac{X}{2t}-\frac{1}{2}\log\ell+C+\psi_d(x,t)\\
 & =\left(1-2\left(\frac{n}{2}-d\right)\right)\log L-\frac{1}{2}\log\ell+C+\psi_d(x,t)
\end{split}
\end{equation}
for $X=2Lt$ and large enough $t$. 
On the other hand, it follows from Lemma~\ref{Lemma:3.1} that 
$$
\lim_{t\to\infty}v_d(2Lte_1,t)=\infty. 
$$
Then, applying by the intermediate value theorem, 
for any large enough $T>0$, 
we find $t_T\in[T,\infty)$ and $j=1,2,\dots$ such that 
\begin{equation}
\label{eq:4.27}
r_j=v_d(2Lt_T,t_T).
\end{equation}
Taking large enough $\ell$ if necessary, 
by \eqref{eq:4.24} and \eqref{eq:4.26} we see that 
\begin{equation*}
\begin{split}
 & \left(\frac{n}{2}-d\right)\log(4\pi t_T)-\log w(2Lt_Te_1,t_T)-\nu_F(v_d(2Lt_Te_1,t_T))\\
 & \le\left(1-2\left(\frac{n}{2}-d\right)\right)\log L-\frac{1}{2}\log\ell+C+\psi_d(2Lt_Te_1,t_T)<-\frac{1}{4}\log\ell.
\end{split}
\end{equation*}
Then, by \eqref{eq:4.21} we obtain 
$$
J_F[v_d](2Lt_Te_1,t;e_1)
\le -\frac{1}{4}\log\ell\cdot \frac{1}{2t_T}(1+o_T(1))<0
$$
for large enough $T$. 
This together with \eqref{eq:4.2}, \eqref{eq:4.25}, and \eqref{eq:4.27} implies that 
$U_d[\phi_\ell](t)$ does not possess the eventual $F$-concavity property. 
Thus Theorem~\ref{Theorem:1.6}~(3) follows. 
$\Box$\vspace{5pt}
\newline
{\bf Proof of Theorem~\ref{Theorem:1.6}~(4).}
Let $d=n/2$ and $\phi\in{\mathcal L}_A$. 
Similarly to the proof of Theorem~\ref{Theorem:1.6}, it suffices to consider the case of $a=1$. 
\vspace{3pt}
\newline
\underline{Step 1}.
Assume that \eqref{eq:1.10} holds. 
It follows from Lemma~\ref{Lemma:2.7} that \eqref{eq:4.5} holds. 
Proposition~\ref{Proposition:4.2} implies that $U_d[\phi]$ possesses the eventual $F$-concavity property. 
\vspace{3pt}
\newline
\underline{Step 2}.
Assume that \eqref{eq:1.10} does not hold. Then $\kappa_F^*\ge 1/2$. 
By Lemma~\ref{Lemma:2.7} we see that 
$\limsup_{r\to\infty}\nu_F(r)>-\infty$. 
Let $\phi_\ell$ be as in the proof of Step~2 of the proof of Theorem~\ref{Theorem:1.6}~(2) with $\ell\in(0,1)$. 
Then $\phi\in{\mathcal L}_A$ and $M_{\phi_\ell}=1/2<1$. 
By Theorem~\ref{Theorem:1.3}~(2) it suffices to consider the case of $\kappa_F\le 1/2$ in $(0,\infty)$. 
Then $\kappa_F^*=1/2$, and 
applying the same argument as in Step 2 in the proof of Theorem~\ref{Theorem:1.6}~(3), 
we see that hat $U_d[\phi_\ell]$ does not possess the eventual $F$-concavity property. 
Thus Theorem~\ref{Theorem:1.6}~(4) follows. 
The proof of Theorem~\ref{Theorem:1.6} is complete.
$\Box$
\vspace{5pt}

Finally, as an application of Theorem~\ref{Theorem:1.6}, 
we prove Corollaries~\ref{Corollary:1.2} and \ref{Corollary:1.3}.
\vspace{5pt}
\newline
{\bf Proof of Corollary~\ref{Corollary:1.2}}.
If the $F$-concavity property is preserved by the heat flow, then
Theorem~\ref{Theorem:1.6} implies that $\limsup_{r\to\infty}\sigma_F(r)<\infty$. 
Then assertion~(1) easily follows from Theorem~\ref{Theorem:1.6}~(2) and (3). 
Since the $1/2$-log-concavity and the hot-concavity properties are preserved by the heat flow in ${\mathbb R}^n$, 
assertion~(2) follows from assertion~(1). 

On the other hand, in the case of $1/2$-log-concavity, 
it follows from Remark~\ref{Remark:1.1}~(2) that $\kappa_F=1/2$ in $(0,\infty)$, that is, $\sigma_F=0$ in $(0,\infty)$. 
Then Theorem~\ref{Theorem:1.6}~(4) implies that 
$U_d[\phi]$ does not possess the eventual $1/2$-log-concavity property for some $\phi\in{\mathcal L}_A$ with $M_\phi<1$. 
Since hot-concavity is stronger than $1/2$-log-concavity in ${\mathcal A}([0,1))$, 
we also see that $U_d[\phi]$ does not possess the eventual hot-concavity property. 
Thus assertion~(3) holds, and the proof of Corollary~\ref{Corollary:1.2} is complete.
$\Box$
\vspace{5pt}
\newline
{\bf Proof of Corollary~\ref{Corollary:1.3}}. 
The eventual $F$-concavity properties which the heat flow develops in ${\mathcal L}_A$  are characterized by Theorem~\ref{Theorem:1.6} with $d=0$. 
It remains to obtain an example of $F$ with properties~(1) and (2). 

Let $\kappa_H$ be the function $\kappa_F$ corresponding to hot concavity
(see (E3)). Since hot concavity is stronger than $1/2$-log-concavity in ${\mathcal A}([0,1))$, 
by Theorem~\ref{Theorem:1.6}~(1), Lemma~\ref{Lemma:2.6}, and Remark~\ref{Remark:1.1}~(3) we see that
\begin{equation}
\label{eq:4.28}
\limsup_{r\to\infty}\, r\left(\kappa_H(r)-\frac{1}{2}\right)<\infty,
\qquad
\kappa_H(r)\ge\frac{1}{2}\quad\mbox{for}\quad r\in(0,\infty).
\end{equation}
Let $k>0$, and set 
$$
f(r):=r^{-1}\max\left\{k(\log r),r\left(\kappa_H(r)-\frac{1}{2}\right)\right\}\quad\mbox{for}\quad r\in(0,\infty).
$$
By \eqref{eq:4.28} we see that $f\ge 0$ in $(0,\infty)$ and  
$f(r)=kr^{-1}(\log r)$ for large enough $r$.
Then we can define an admissible function $F$ on $[0,1)$ by 
$$
F(\tau):=-\int_1^{-\log\tau} s^{-\frac{1}{2}}\exp\left(\int^\infty_s \xi^{-1}f(\xi)\,d\xi\right)\,ds\quad\mbox{for}\quad \tau\in(0,1),
\qquad
F(0)=-\infty.
$$
It follows that 
\begin{align*}
{\mathcal F}(r) & =F(e^{-r})=-\int_1^r s^{-\frac{1}{2}}\exp\left(\int^\infty_s \xi^{-1}f(\xi)\,d\xi\right)\,ds,\\
{\mathcal F}'(r) & =-r^{-\frac{1}{2}}\exp\left(\int^\infty_r \xi^{-1}f(\xi)\,d\xi\right)<0,\\
{\mathcal F}''(r) & =\left(\frac{1}{2}r^{-\frac{3}{2}}+r^{-\frac{3}{2}}f(r)\right)\exp\left(\int^\infty_r \xi^{-1}f(\xi)\,d\xi\right),
\end{align*}
for $r\in(0,\infty)$. Then we see that $F'>0$ in $(0,1)$ and 
$$
\kappa_F(r)=\frac{1}{2}+f(r)\ge \kappa_H(r)\quad\mbox{for}\quad r\in(0,\infty),
\quad
\kappa_F(r)>\kappa_H(r)\quad\mbox{for large enough $r$.}
$$
This together with Lemma~\ref{Lemma:2.6} implies that $F$-concavity is strictly stronger than hot-concavity in ${\mathcal A}([0,1))$. 
Furthermore, by \eqref{eq:4.28} we have 
$$
\sigma_F(r)=rf(r)=k\log r\quad\mbox{for large enough $r>0$}.
$$
Setting $k\in(0,1/4)$, by Theorem~\ref{Theorem:1.6}~(2), (3) 
we see that $e^{t\Delta}\phi$ possesses the eventual $F$-concavity property for $\phi\in{\mathcal L}_A$.
Thus Corollary~\ref{Corollary:1.3} follows.
$\Box$
\medskip

\noindent
{\bf Acknowledgements.} 
The author of this paper thanks Paolo Salani and Asuka Takatsu for fruitful discussions. 
He was supported in part by JSPS KAKENHI Grant Number 19H05599. 
\medskip 

\noindent
{\bf Data availability.}
No data was generated or analyzed as part of the writing of this paper.
\medskip 

\noindent{\bf Declarations.}\\
{\bf Conflict of interest}: 
The author of this paper states that there is no conflict of interest.


\begin{bibdiv}
\begin{biblist}
\bib{AV}{article}{
   author={Aronson, D. G.},
   author={V\'{a}zquez, J. L.},
   title={Eventual $C^\infty$-regularity and concavity for flows in
   one-dimensional porous media},
   journal={Arch. Rational Mech. Anal.},
   volume={99},
   date={1987},
   pages={329--348},
}
\bib{BL}{article}{
 author={Brascamp, Herm Jan},
 author={Lieb, Elliott H.},
 title={On extensions of the Brunn-Minkowski and Pr\'{e}kopa-Leindler
 theorems, including inequalities for log concave functions, and with an
 application to the diffusion equation},
 journal={J. Functional Analysis},
 volume={22},
 date={1976},
 pages={366--389},
}
\bib{BV}{article}{
   author={B\'{e}nilan, Philippe},
   author={V\'{a}zquez, Juan Luis},
   title={Concavity of solutions of the porous medium equation},
   journal={Trans. Amer. Math. Soc.},
   volume={299},
   date={1987},
   pages={81--93},
}
\bib{CW1}{article}{
   author={Chau, Albert},
   author={Weinkove, Ben},
   title={Counterexamples to quasiconcavity for the heat equation},
   journal={Int. Math. Res. Not. IMRN},
   date={2020},
   pages={8564--8579},
}
\bib{CW2}{article}{
   author={Chau, Albert},
   author={Weinkove, Ben},
   title={The Stefan problem and concavity},
   journal={Calc. Var. Partial Differential Equations},
   volume={60},
   date={2021},
   pages={Paper No. 176, 13},
}
\bib{CW3}{article}{
   author={Chau, Albert},
   author={Weinkove, Ben},
   title={Non-preservation of $\alpha$-concavity for the porous medium equation},
    journal={preprint (arXiv:2011.03063)}
}
\bib{GK}{article}{
   author={Greco, Antonio},
   author={Kawohl, Bernd},
   title={Log-concavity in some parabolic problems},
   journal={Electron. J. Differential Equations},
   date={1999},
   pages={No. 19, 12},
}
\bib{INS}{article}{
   author={Ishige, Kazuhiro},
   author={Nakagawa, Kazushige},
   author={Salani, Paolo},
   title={Spatial concavity of solutions to parabolic systems},
   journal={Ann. Sc. Norm. Super. Pisa Cl. Sci.},
   volume={20},
   date={2020},
   pages={291--313},
}
\bib{IS01}{article}{
 author={Ishige, Kazuhiro},
 author={Salani, Paolo},
 title={Is quasi-concavity preserved by heat flow?},
 journal={Arch. Math. (Basel)},
 volume={90},
 date={2008},
 pages={450--460},
}
\bib{IS02}{article}{
 author={Ishige, Kazuhiro},
 author={Salani, Paolo},
 title={Convexity breaking of the free boundary for porous medium
 equations},
 journal={Interfaces Free Bound.},
 volume={12},
 date={2010},
 pages={75--84},
}
\bib{HV}{article}{
   author={Huang, Yong},
   author={V\'{a}zquez, Juan L.},
   title={Large-time geometrical properties of solutions of the Barenblatt equation of elasto-plastic filtration},
   journal={J. Differential Equations},
   volume={252},
   date={2012},
   pages={4229--4242},
}
\bib{IST01}{article}{
 author={Ishige, Kazuhiro},
 author={Salani, Paolo},
 author={Takatsu, Asuka},
 title={To logconcavity and beyond},
 journal={Commun. Contemp. Math.},
 volume={22},
 date={2020},
 pages={1950009, 17},
}
\bib{IST03}{article}{
   author={Ishige, Kazuhiro},
   author={Salani, Paolo},
   author={Takatsu, Asuka},
   title={New characterizations of log-concavity via Dirichlet heat flow},
   journal={Ann. Mat. Pura Appl.},
   volume={201},
   date={2022},
   pages={1531--1552},
}
\bib{IST05}{article}{
 author={Ishige, Kazuhiro},
 author={Salani, Paolo},
 author={Takatsu, Asuka},
 title={Characterization of $F$-concavity preserved by the Dirichlet heat flow},
 journal={preprint (arXiv:2207.13449)},
}
\bib{IST06}{article}{
 author={Ishige, Kazuhiro},
 author={Salani, Paolo},
 author={Takatsu, Asuka},
 title={Hierarchy of deformations in concavity},
 journal={to appear in Inf. Geom.},
}
\bib{Keady}{article}{
   author={Keady, G.},
   title={The persistence of logconcavity for positive solutions of the
   one-dimensional heat equation},
   journal={J. Austral. Math. Soc. Ser. A},
   volume={48},
   date={1990},
   pages={246--263},
}
\bib{KL}{article}{
   author={Kim, Sunghoon},
   author={Lee, Ki-Ahm},
   title={Local continuity and asymptotic behaviour of degenerate parabolic systems},
   journal={Nonlinear Anal.},
   volume={192},
   date={2020},
   pages={111702, 32},
}
\bib{Ken}{article}{
 author={Kennington, Alan U.},
 title={Power concavity and boundary value problems},
 journal={Indiana Univ. Math. J.},
 volume={34},
 date={1985},
 pages={687--704},
}
\bib{Ken02}{article}{
   author={Kennington, Alan U.},
   title={Convexity of level curves for an initial value problem},
   journal={J. Math. Anal. Appl.},
   volume={133},
   date={1988},
   pages={324--330},
}
\bib{Kol}{article}{
   author={Kolesnikov, Alexander V.},
   title={On diffusion semigroups preserving the log-concavity},
   journal={J. Funct. Anal.},
   volume={186},
   date={2001},
   pages={196--205},
}
\bib{Kor}{article}{
 author={Korevaar, Nicholas J.},
 title={Convex solutions to nonlinear elliptic and parabolic boundary
 value problems},
 journal={Indiana Univ. Math. J.},
 volume={32},
 date={1983},
 pages={603--614},
}
\bib{L}{article}{
   author={Lee, Ki-Ahm},
   title={Power concavity on nonlinear parabolic flows},
   journal={Comm. Pure Appl. Math.},
   volume={58},
   date={2005},
   pages={1529--1543},
}
\bib{LPV}{article}{
   author={Lee, Ki-Ahm},
   author={Petrosyan, Arshak},
   author={V\'{a}zquez, Juan Luis},
   title={Large-time geometric properties of solutions of the evolution
   $p$-Laplacian equation},
   journal={J. Differential Equations},
   volume={229},
   date={2006},
   pages={389--411},
}
\bib{LV}{article}{
   author={Lee, Ki-Ahm},
   author={V\'{a}zquez, J. L.},
   title={Geometrical properties of solutions of the porous medium equation
   for large times},
   journal={Indiana Univ. Math. J.},
   volume={52},
   date={2003},
   pages={991--1016},
}
\end{biblist}
\end{bibdiv}
\end{document}